\documentclass[11pt]{article}
\usepackage{latexsym,amssymb,enumerate}
\usepackage{amsmath}
\usepackage{color}

\newtheorem{Theorem}{\bf Theorem}[section]
\newtheorem{Lemma}{\bf Lemma}[section]
\newtheorem{Proposition}{\bf Proposition}[section]
\newtheorem{Corollary}{\bf Corollary}[section]
\newtheorem{Remark}{\bf Remark}[section]
\newtheorem{Example}{\bf Example}[section]
\newtheorem{Definition}{\bf Definition}[section]

\newenvironment{theorem}{\begin{Theorem}$\!\!\!$}{\end{Theorem}}
\newenvironment{lemma}{\begin{Lemma}$\!\!\!$}{\end{Lemma}}
\newenvironment{proposition}{\begin{Proposition}$\!\!\!$}{\end{Proposition}}

\newenvironment{remark}{\begin{Remark}$\!\!\!$}{\end{Remark}}

\newenvironment{definition}{\begin{Definition}$\!\!\!$}{\end{Definition}}

\def\exp{\mathrm {exp}}

\let\lt=<
\let\gt=>

\def\R{{\mathbb R}}

\def\N{{\mathbb N}}

\let\eps=\varepsilon

\def\be{\begin{equation}}
\def\ee{\end{equation}}
\def\beq{\begin{eqnarray*}}
\def\eeq{\end{eqnarray*}}
\def\e{{\rm e}}

\numberwithin{equation}{section}

\begin{document}

\title{Asymptotic  Behavior and Decay Estimates of  the Solutions for 
a Nonlinear Parabolic Equation with Exponential Nonlinearity}

\author
{
Giulia FURIOLI  \\
{\small Universit\`a di Bergamo, Italy}\\
\\
Tatsuki KAWAKAMI  \\
{\small Osaka Prefecture University, Japan}\\
\\
 Bernhard RUF\\
{\small Universit\`a di Milano, Italy}\\
\\
Elide TERRANEO\\
{\small Universit\`a di Milano, Italy}
\\
}

\date{}
\maketitle
\begin{abstract} \noindent
We consider a nonlinear parabolic equation with an exponential nonlinearity 
which is critical with respect to the growth of the nonlinearity and the regularity of the initial data.
After showing the equivalence of the notions of weak and mild solutions, we derive
decay estimates and the asymptotic behavior for small global-in-time solutions.
\end{abstract}
 
\par \bigskip
\section{Introduction}
We consider the Cauchy problem for the semilinear heat equation
\begin{equation}
\label{eq:1.1}
\left\{
\begin{array}{ll}
\partial_t u= \Delta u + f(u), \quad
& t>0,\  x\in\R^n,
\vspace{5pt}\\
u(0,x)=\varphi (x) , \quad
&x\in{\mathbb R}^n,
\end{array}
\right. \ee
where $n\ge 1$, $u(t,x):\R^+\times \R^n\rightarrow \R$ is the unknown function,
$\partial_t=\partial/\partial t$, 
$f(u)$ contains the nonlinearity and $\varphi$ is the initial data. 
This paper is concerned with the asymptotic behavior and decay estimates of the solutions of \eqref{eq:1.1}
in certain limiting cases which are critical with respect to the growth of the nonlinearity and
the regularity of the initial data.
\par \medskip 
Before introducing the subject of this paper, let us recall some related known results.
\par \medskip \noindent
\noindent
{\it{The polynomial case.}}
The case of power nonlinearities $f(u)=|u|^{p-1}u$ with $p>1$, that is
 \begin{equation}
 \label{eq:1.2}
\left\{
\begin{array}{ll}
\partial_t u= \Delta u + |u|^{p-1}u, \quad
&t>0,\  x\in\R^n,
\vspace{5pt}\\
u(0,x)=\varphi (x),\quad
&x\in{\mathbb R}^n,
\end{array}
\right.
\end{equation}
has been extensively studied since the pioneering works by Fujita \cite{F} and Weissler \cite{W1,W2}.
It is well-known that the problem \eqref{eq:1.2} satisfies a scale invariance property.
In fact, for $\lambda\in{\mathbb R}_+$,
if $u$ is a solution of \eqref{eq:1.2}, 
then  
\begin{equation}
\label{eq:1.3}
u_\lambda(t,x) := \lambda^{\frac 2{p-1}} u(\lambda^2 t, \lambda x)
\end{equation}
is also a solution of \eqref{eq:1.2} with initial data 
$\varphi_\lambda(x):=\lambda^{2/(p-1)}\varphi(\lambda x)$. 
So, all function spaces invariant with respect to the scaling transformation \eqref{eq:1.3}
play a fundamental role in the study of the Cauchy problem \eqref{eq:1.2}. 
In the framework of Lebesgue spaces, 
we can easily show that the norm of the space $L^q({\mathbb R}^n)$ is invariant with respect to \eqref{eq:1.3}
if and only if $q=q_c$ with $q_c=n (p-1)/2$, and it is well-known that, given the nonlinearity $|u|^{p-1}u$, the Lebesgue space $L^{q_c}({\mathbb R}^n)$ 
plays the role of {\it critical space} for the well-posedness 
of \eqref{eq:1.2} (see e.g. \cite{BC,HW,W1,W2}).  
\par 
Indeed, for any $q\geq q_c$ and $q>1$, or $q>q_c$ and $q\geq 1$ ({\it subcritical case}),
 and for any initial data $\varphi\in L^q({\mathbb R}^n)$
 there exists a local (in time) solution of the Cauchy problem
 \eqref{eq:1.2}. On the other hand, for initial data belonging to $L^q({\mathbb R}^n)$ with $1<q<q_c$ 
 ({\it supercritical case}), 
 Weissler \cite {W1} and Brezis-Cazenave \cite{BC} indicate that for certain $\varphi \in L^q({\mathbb R}^n)$
  there exists no local (in time) solution in any  reasonable  sense.
\par 
We can also state the previous results from a different point of
view.  Given any initial data in the  Lebesgue space $L^q$,
$1<q<+\infty$, the Cauchy problem \eqref{eq:1.2} is well-posed if
and only if the power nonlinearity $p$ is smaller than or equal
to the critical value $ p_c=1+(2q)/n$. Moreover, the {\it critical case},
which is defined equivalently by $q=q_c$ or $p=p_c$, is the only
case for which global (in time) existence can be established
for small initial data.
\par \medskip
The same polynomial nonlinearity has also been considered for the
Schr\"{o}dinger equation
 \begin{equation}
 \label{eq:1.4}
\left\{
\begin{array}{ll}
i\partial_t u+ \Delta u = |u|^{p-1}u, \quad &t>0,\  x\in\R^n,
\vspace{5pt}\\
u(0,x)=\varphi (x),\quad &x\in{\mathbb R}^n.
\end{array}
\right.
\end{equation}
Here $u:\mathbb R^+\times \mathbb R^n\to \mathbb C$ and 
$\varphi \in H^s(\mathbb R^n)$ with $0\leq s<n/2$.  
Also for problem \eqref{eq:1.4}
the scaling invariance \eqref{eq:1.3} holds,
and the scaling argument indicates the critical Sobolev space and exponent for the nonlinearity $|u|^{p-1}u$:
for any fixed $p>1$,  the Sobolev space $H^{s_c}({\mathbb R}^n)$ with $s_c=n/2-2/(p-1)$ plays the role of {\it critical space}
for the well-posedness of \eqref{eq:1.4} (see e.g. \cite{CW}).
 Indeed, for initial data $\varphi\in H^s({\mathbb R}^n)$ with $s \ge s_c$ ({\it subcritical case})   
 Cazenave and Weissler \cite{CW} proved that there exists a local (in time) solution in $H^s(\mathbb R^n)$. 
\par
Equivalently, for any initial data in the Sobolev space $H^s({\mathbb R}^n)$, $0\leq s<n/2$, 
local existence for the Cauchy problem \eqref{eq:1.4}  holds true under the condition
$p\leq \tilde p_c=1+ {4}/(n-2s)$, namely for any power nonlinearity
smaller than or equal to the critical nonlinearity $\tilde p_c$.

Moreover, the {\it critical case}, which is defined equivalently by
$s=s_c$ or $p=\tilde p_c$, is the only one for which global (in time)
existence can be established for small initial data.
\par 
The case $s_c =  n/ 2$ is a {\it limiting case}, since for $s >n/ 2$ the Sobolev space $H^s(\mathbb R^n)$ embeds into $L^\infty(\mathbb R^n)$, and no specific behavior for the nonlinearity is required for the existence of local-in-time solutions of  \eqref{eq:1.4} in the $H^s$ framework.
\vspace{8pt}
 
\noindent
{\it{The limiting case.}} Let us now consider the {\it limiting critical case} $H^{n/2}(\mathbb R^n)$; then
\par \noindent \smallskip
(i) any power nonlinearity is subcritical, since $H^{n/2}(\mathbb R^n)$ embeds into any $L^p(\mathbb R^n)$ space,
\par \noindent
(ii) $H^{n/2}(\mathbb R^n)$ does not embed into $L^{\infty}(\mathbb R^n)$, 
and by Trudinger's inequality \cite{T} one knows that $H^{n/2}(\mathbb R^n)$ 
embeds into the Orlicz space $L_\varphi (\mathbb R^n)$, 
with growth function ($N$-function) $\varphi(s) = e^{s^2} - 1$.
\par \smallskip
For this limiting case,
Nakamura and Ozawa \cite{NO} considered the Cauchy problem
 \begin{equation}
 \label{eq:1.5}
\left\{
\begin{array}{ll}
i\partial_t u+ \Delta u = f(u), \quad & t>0,\  x\in\R^n,
\vspace{5pt}\\
u(0)=\varphi\in H^{\frac{n}{2}}({\mathbb R}^n), &\,
\end{array}
\right.
\end{equation}
with exponential nonlinearity of asymptotic growth $f(u)\sim e^{u^2}$ and  
with a vanishing behavior at the origin. 
They proved the existence of global-in-time solutions of \eqref{eq:1.5}
for initial data with small norm in $H^{n/2}(\mathbb R^n)$.
In view of the Trudinger inequality
the growth rate of $f(u)$ at infinity seems to be optimal for initial data in $H^{n/2}(\mathbb R^n)$.
\par \medskip
It is well-known that  there is a correspondence
between the $H^s({\mathbb R}^n)$ theory for the nonlinear Schr\"{o}dinger equation \eqref{eq:1.4} 
and the $L^q({\mathbb R}^n)$ theory for the semilinear heat equation \eqref{eq:1.2}
(see e.g. \cite{CW}).
Indeed,  since the Sobolev space $H^s(\mathbb R^n)$ embeds into $ L^q(\mathbb R^n)$ with $q=(2n)/(n-2s)$,
for $0\leq s<n/2$,
it holds that
$p_c=1+(2q)/n=1+ 4/(n-2s)=\tilde p_c$.
This means that
the critical nonlinearity $ |u|^{p_c-1}u$ associated to \eqref{eq:1.2} in the $L^q({\mathbb R}^n)$ framework
corresponds to the critical one associated to \eqref{eq:1.4} in the $H^s({\mathbb R}^n)$ framework.

As a natural analogy to the results of \cite{NO},
the third and fourth authors of this paper \cite{RT} and Ioku \cite{I,I-t} considered 
the Cauchy problem \eqref{eq:1.1} with a nonlinearity of the form
\be
\label{eq:1.6}
f(u)= |u|^{\frac{4}{n}}u\,{e}^{u^2}
\ee 
and initial data $\varphi$ belonging to the Orlicz space $\exp{L^2}$ defined as
$$
\exp{L^2}(\mathbb R^n)
:=\left\{ u\in L^1_{loc}(\mathbb R^n),\ \exists \lambda>0 :
\int_{\mathbb R^n}\left({\exp}\left(\frac{|u(x)|}{\lambda}\right)^2-1\right)dx<\infty\right\}
$$
(see also Definition~\ref{Definition:3.1}). Clearly, no scaling invariance holds true for equation 
\eqref{eq:1.1} with such a nonlinear term. Also, the growth $e^{u^2}$ at infinity of the nonlinearity $f(u)$ 
seems to be optimal in the framework of the Orlicz space $\exp L^2$. 
In fact, if $f(u)\sim {e}^{|u|^r}$ with $r>2$,
then there exist positive initial data $\varphi\in \exp L^2$ (even with very small norm)
such that there exists no nonnegative classical local-in-time solution of \eqref{eq:1.1} (see \cite{IRT}).

For the Cauchy problem \eqref{eq:1.1} with \eqref{eq:1.6},
the authors of the papers \cite{I, I-t, RT} considered the corresponding integral equation
\begin{equation}
\label{eq:1.7}
u(t)= e^{t\Delta} \varphi + \int_{0}^te^{(t-s)\Delta} f(u(s)) ds.
\end{equation}
As we will recall in Definition \ref{mild}, $u$ is called a {\it mild} solution of \eqref{eq:1.1}  if $u$ is a solution of the integral equation \eqref{eq:1.7}
and $u(t) \underset {t\to 0}{\longrightarrow}\varphi$ in the weak$^*$ topology (see \eqref{eq:2.1}).
Under this notion of solution one has
\begin{proposition}\ {\rm \cite{I,I-t,RT}}
\label{Proposition:1.1}
Let $n\ge1$ and $\varphi \in \exp{L^2}$.
Suppose that $f$ satisfies \eqref{eq:1.6}.
Then there exists  $\varepsilon=\eps(n) >0$ such that, 
if $\|\varphi\|_{\exp{L^2}} <\eps$, 
then there exists a mild solution $u$ of \eqref{eq:1.1}
satisfying
$$
u\in L^\infty(0,\infty\,;\, \exp{L^2})
$$
and
\begin{equation}
\label{eq:1.8}
\sup _{t>0} \|u(t)\|_{\exp{L^2}} \leq 2\|\varphi\|_{\exp{L^2}}.
\end{equation}
\end{proposition}

By these results one obtains the small data global existence of {\it mild} solutions of \eqref{eq:1.1}.
However, as far as we know, there are no results which treat the correspondence between mild solutions and other notions of weak solutions, and which prove decay estimates and the asymptotic behavior of the solutions of \eqref{eq:1.1}.
\par \medskip
In this paper, under a smallness assumption for the solution,
we prove the equivalence between {\it mild} solutions and  {\it weak} solutions of \eqref{eq:1.1}.
Furthermore, 
under condition \eqref{eq:1.8},
we obtain decay estimates for the solutions in the following two cases
\par
$\varphi\in \exp L^2$ only ({\it singular case}), and
\par 
$\varphi \in \exp L^2\cap L^p(\mathbb R^n)$ with $p\in [1,2)$ ({\it regular case}).
\par \noindent
In particular, for the regular case $p=1$, 
we show that global-in-time solutions with some suitable decay estimates behave asymptotically
like suitable multiples of the heat kernel.
\vspace{5pt}

The paper is organized as follows.
In Section~2 
we state the main results.
In Section~3
we give some preliminaries. 
Sections~4--7
are devoted to the proof of the theorems,
respectively.
In Section~8
we extend our theorems to fractional diffusion equations with general initial data.
\vspace{8pt}

Before closing this section
we give some notations used in this paper.
For any $p\in[1,\infty]$, let $\|\cdot\|_{L^p}$ be the usual norm of $L^p:=L^p({\mathbb R}^n)$; we 
write $\mathbb N_0:=\mathbb N\cup\{0\}$;
$\mathcal{D}(\Omega)$ denotes the space of $C^\infty$- functions with compact support in $\Omega$,
and $\mathcal{D}'(\Omega)$ the topological dual of $\mathcal{D}(\Omega)$.
We denote by $\mathcal{S}'(\Omega)$ the space of tempered distributions.
Throughout the present paper,
$C$ will denote a generic positive constant which may have different values also within the same line.

\section{Main results}
In this section we state the main results of this paper.
Before presenting the main theorems
we introduce the definition of {\it weak} and  {\it mild} solutions.
\begin{definition} {\rm(Weak solution)}
\label{Definition:2.1}
Let  $\varphi\in \exp L^2$, $T\in (0,\infty]$ and $u\in L^{\infty}(0,T; \exp L^2)$.
We call $u$ a weak solution of the Cauchy problem \eqref{eq:1.1} in $(0,T)\times{\mathbb R}^n$ 
if $u$ satisfies
\[
\partial_t u= \Delta u + f(u)
\]
in the sense of distributions ${\cal D}'((0,T)\times \R^n)$
and $u(t) \underset {t\to 0}{\longrightarrow} \varphi$ in weak$^*$
topology.
\end{definition}
\begin{definition} {\rm(Mild solution)}
\label{mild}
Let  $\varphi\in \exp L^2$, $T\in (0,\infty]$ and $u\in L^\infty(0,T; \exp L^2)$.
We call $u$ a mild solution of the Cauchy problem \eqref{eq:1.1} in $(0,T)\times{\mathbb R}^n$
if $u\in C((0,T); \exp L^2)$, $u$ is a solution of the integral equation 
\[
u(t)= e^{t\Delta} \varphi + \int_{0}^te^{(t-s)\Delta} f(u(s)) ds
\]
and $u(t) \underset {t\to 0}{\longrightarrow}\varphi$ in the weak$^*$ topology.
\end{definition}
We recall that $u(t) \underset {t\to 0}{\longrightarrow}
\varphi$ in weak$^*$ topology if and only if
\begin{equation}\label{eq:2.1} \lim_{t\to 0}
\int_{\mathbb{R}^n}\bigl(u(t,x)\psi(x) -\varphi(x)\psi(x)\bigr) dx=0
\end{equation}
{for  every} $ \psi$ belonging to
 the predual space of
$\exp L^2$ (see  Section 3).
\vspace{5pt}

Now we are ready to state the main results of this paper.
First we show the equivalence between 
 small {\it mild} solutions and small {\it weak} solutions of \eqref{eq:1.1}.
\begin{proposition}
\label{Theorem:2.1}
Let $\varphi\in\exp L^2$, $T\in(0,\infty]$ and $u\in L^{\infty}(0,T; \exp L^2)$.
Then there exists a positive constant $\varepsilon>0$ such that,
if $\sup_{0<t<T}\|u(t)\|_{\exp L^2}\leq \varepsilon$,
then the following statements are equivalent:
\begin{itemize}
\item[{\rm(i)}] $u$ is a weak solution of \eqref{eq:1.1} in $(0,T)\times{\mathbb R}^n$;
\item[{\rm(ii)}] $u$ is a mild solution of \eqref{eq:1.1} in $(0,T)\times{\mathbb R}^n$.
\end{itemize}
\end{proposition}
Here and in the rest of the paper, 
for simplicity,
we call $u$ a solution of \eqref{eq:1.1}
if $u$ is a mild solution of \eqref{eq:1.1} in $(0,\infty)\times{\mathbb R}^n$.
\vspace{5pt}

Next we prove uniqueness of  small solutions in $L^\infty((0,T); \exp L^2)$ for all dimensions $n\in\N$.
\begin{proposition}
\label{uniqueness}
Let $n\geq 1$.  There is $\eps >0$ such that if  $\varphi\in\exp L^2(\R^n)$ with $\|\varphi\|_{\exp L^2}\leq \eps$ and   $u$, $v$ are two  solutions of \eqref{eq:1.1} satisfying 
\be\label{smallness}
\sup_{t<T} \|u(t)\|_{\exp L^2} \leq 2 \|\varphi\|_{\exp L^2}, \quad \sup_{t<T'} \|u(t)\|_{\exp L^2} \leq 2 \|\varphi\|_{\exp L^2},
\ee
then $u(t)=v(t)$ on $[0, \min (T, T'))$.
\end{proposition}
The following result provides a decay estimate of the small solution of \eqref{eq:1.1} 
for the singular case, that is, $\varphi\in\exp L^2$ only.
\begin{theorem} 
\label{Theorem:2.2}
Let $n\geq 1$ and $\varphi\in\exp L^2$ with $\varphi\ge0$.
Assume that there exists a unique positive solution $u$ of \eqref{eq:1.1} satisfying \eqref{eq:1.8}. 
Then there exist constants $\eps = \eps (n) >0$ and $C=C(n)>0$
such that, 
if $\|\varphi\|_{\exp L^2} <\eps$, 
then the solution $u$ satisfies
\be
\label{eq:2.2} 
\|u(t)\|_{L^q} \leq C t^{-\frac n2 \left( \frac 12-\frac 1q\right )}\|\varphi\|_{\exp L^2},\quad t>0,
\ee
for all $q\in[2,\infty]$.
\end{theorem}
\begin{remark}
{\rm(i)}
By Propositions~$\ref{Proposition:1.1}$ and ~$\ref{uniqueness}$,
if $\|\varphi\|_{\exp L^2}$ is small enough,
then we can show that the assumption for the solution $u$ is not empty. 
\vspace{3pt}
\newline
{\rm(ii)}
We obtained the same decay estimate as the solution of the linear heat equation with initial data in $L^2$. 
\vspace{3pt}
\newline
{\rm(iii)}
For the lower dimensional case $1\le n\le4$, 
in the proof of Proposition~$\ref{Proposition:1.1}$,
Ioku {\rm\cite{I,I-t}} already obtained the decay estimate \eqref{eq:2.2}
for $1+4/n<q<2+8/n$.
\vspace{3pt}
\newline
{\rm(iv)}
Due to the embedding  $\exp L^2 \subset L^q$ for $2\leq q<\infty$ 
{\rm(}see Lemma~$\ref{Lemma:3.2}${\rm)}, 
if $\varphi \in \exp{L^2}$ with
$\|\varphi\|_{\exp{L^2}} <\eps$ as in Proposition~$\ref{Proposition:1.1}$, 
then $u\in L^\infty(0,\infty\,;\, L^q)$ for all $2\leq q <\infty$.
\end{remark}

Next we consider the regular case, namely $\varphi\in\exp L^2 \cap L^p$, $p\in [1,2)$.
The nonlinearity \eqref{eq:1.6} satisfies  $f(x)\sim x^{1+4/n}$ for $x\to 0^+$. 
So, if $u\in L^\infty(0,\infty; L^q)$ for $q\geq 2$,
then the nonlinear term $f(u)$ belongs to  $L^p$ for $p \geq (2n)/(n+4)$.
In the lower dimensional case $1\leq n\leq 4$, this means that $f(u)\in L^p$ for all $p\geq 1$
but for the higher dimensional case $n\geq 5$ this implies a true constraint.
This is the reason why in the next theorems we have to introduce some parameters 
$p_*$, $p_1$ (and $p_2$ in Lemma \ref{Lemma:6.1}) 
which are meaningful only for the higher dimensional case. 
More specifically, in the higher dimensional case, 
we can prove two kinds of results about the decay estimate of solutions of \eqref{eq:1.1}. 
In the first one (Theorem \ref{Theorem:2.3}), 
we only assume a control of the $\exp L^2$ norm of the initial data,
allowing the $L^p$-norm of the same data to be large (see Remark~\ref{Remark:6.3}). 
This mild assumption on the initial data entails a decay estimate of $\|u(t)\|_{L^q}$ 
only for $q\geq p_*>p$.
In the second result (Theorem \ref{Theorem:2.4}), 
under a stronger assumption, that is, 
a smallness assumption for both the $\exp L^2$ and the $L^p$ norm of the initial data, 
we can prove a better decay estimate on $\|u(t)\|_{L^q}$ for all $q\geq p$.
\begin{theorem} 
\label{Theorem:2.3}
Assume the same conditions as in Theorem~$\ref{Theorem:2.2}$.
Furthermore, assume $\varphi\in L^p$ for some $p\in[1,2)$.
Put
\begin{equation}
\label{eq:2.3}
p_*:=\max\left\{p,\frac{2n}{n+4}\right\}.
\end{equation}
Then there exist positive constants $\eps =\eps (n) $, $C=C(n)$
and a positive function $F=F(n, p_*, \|\varphi\|_{p_*})$ such that, if
\be\label{eq:2.4}
\|\varphi\|_{\exp L^2} < \min \left(\eps, F(n,p_*, \|\varphi\|_{L^{p_*}})\right)
\ee
 then the solution $u$ satisfies
\be\label{eq:2.5}
\|u(t)\|_{L^q} \leq C t^{-\frac n2 \left( \frac 1{p_*} -\frac 1q\right )}\|\varphi\|_{\exp L^2\cap L^{p_*}},\quad t>0,
\ee
for all $q\in[p_*,\infty]$.
In particular, if $p\in (p_1,2)$, then
\be\label{eq:2.6}
\|u(t)\|_{L^q} = o\left(t^{-\frac n2 \left( \frac 1{p_*} -\frac 1q\right )}\right ),\quad t\to \infty,
\ee
where
\begin{equation}
\label{eq:2.7}
p_1:=\max\left\{1,\frac{2n}{n+4}\right\}.
\end{equation}
\end{theorem}
\begin{theorem} 
\label{Theorem:2.4} 
Assume the same conditions as in Theorem~$\ref{Theorem:2.3}$. 
Then there exists a positive constant $\delta=\delta(n)$ such that, if 
$$
\max\{\|\varphi\|_{\exp L^2}, \|\varphi\|_{L^p}\}<\delta, 
$$
then \eqref{eq:2.5} with $p_*=p$ holds for all $q\in[p,\infty]$. 
In particular, for all $q\in[p,\infty)$,
$$
\|u(t)\|_{L^q} \leq C(1+t)^{-\frac n2 \left( \frac 1p -\frac 1q\right )}\|\varphi\|_{\exp L^2\cap L^p},\quad t>0.
$$
Furthermore, if $p\in (1,2)$, then \eqref{eq:2.6} with $p_*=p$ holds.
\end{theorem}

Finally we address the question of the asymptotic behavior of solutions of \eqref{eq:1.1}
when $\varphi \in\exp L^2 \cap L^1$.
We show that global-in-time solutions with suitable decay properties behave asymptotically
like suitable multiples of the heat kernel.

\begin{theorem} 
\label{Theorem:2.5} 
Let $n\geq 1$, $\varphi \ge 0$ and $\varphi \in\exp L^2 \cap L^1$. 
Assume that $\|\varphi\|_{\exp L^2}$ is small enough.
Furthermore, suppose that
\begin{enumerate}[a)]
\item for $n\geq 1$,
\be
\label{eq:2.8}
\sup_{t>0}\,t^{\frac n2\left(1-\frac 1q\right )}\|u(t)\|_{L^q}<\infty,\qquad q\in[1,\infty];
\ee
\item for $n \geq 5$, assume moreover that there is $T^*>0$ such that
\be
\label{eq:2.9}
\sup_{0<t<T^*}\|u(t)\|_{L^{\frac 4n+1}}<\infty.
\ee
\end{enumerate}
Then there exists the limit
\[
 \lim_{t\to \infty} \int_{\R^n}u(t,x) \, dx
 = \int_{\R^n} \varphi(x) dx + \int_{0}^{\infty} \int_{\R^n} f(u(t,x))\, dx dt :=m_*
\]
such that
 \be
 \label{eq:2.10}
 \lim_{t\to \infty} t^{\frac n2 \left(1 -\frac 1q\right )}\|u(t)-m_* g(t)\|_{L^q}=0,\qquad q\in[1,\infty],
 \ee
 where $g(t,x)=G_2(t+1,x)$ is the heat kernel, that is,
 $$
 G_2(t,x)= (4\pi t)^{-\frac n2} e^{-\frac {|x|^2}{4t}}.
$$
\end{theorem}
\begin{remark}
\label{Remark:2.3}
{\rm(i)} Conditions \eqref{eq:2.8} and \eqref{eq:2.9} are fulfilled for example 
under the hypotheses of Theorems~$\ref{Theorem:2.3}$ and $\ref{Theorem:2.4}$.
\vspace{3pt}
\newline
{\rm(ii)}
The large time behavior and the decay properties of the solution of \eqref{eq:1.2}
have  been widely studied 
{\rm(}see, i.e., {\rm\cite{EK, IKK1,IKK2,Ka,LN,S,Wa}} and the references therein{\rm)}.
In order to treat these topics for \eqref{eq:1.2},
the scale transformation \eqref{eq:1.3} plays an important role.
However, unfortunately, for the problem \eqref{eq:1.1} with \eqref{eq:1.6},
we don't have such a scaling property.
Instead of \eqref{eq:1.3}, applying the embedding $\exp L^2 \subset L^q$ for $2\leq q<\infty$
and exploiting the uniqueness of the solution,
we overcome this difficulty, and this is one of the novelties of this paper.
\end{remark}
%

In the last section of the paper, 
we consider the Cauchy problem associated to the fractional diffusion equation
\begin{equation}
\label{eq:2.11}
 \left\{
\begin{array}{ll}
\partial_t u+{\cal L}_{\theta}u=|u|^{\frac {r\theta}n}u\,{e}^{u^r}, \quad
& t>0,\  x\in\R^n,
\vspace{5pt}\\
u(0,x)=\varphi (x) , \quad &x\in{\mathbb R}^n,
\end{array}
\right.
\end{equation}
where $0<\theta\le2$, $r>1$ and $\varphi\in \exp L^r$
(we will recall the definition of the Orlicz spaces $\exp L^r$ in Definition~\ref{Definition:3.1}).
Here the operator ${\cal L}_{\theta}:=(-\Delta)^{\theta/2}$ is the fractional Laplacian
defined by the Fourier transform $\mathcal{F}$ as
\begin{equation}
\label{eq:2.12}
{\cal L}_{\theta}\phi:= \mathcal{F}^{-1}[|\xi|^{\theta}\mathcal{F}[\phi]].
\end{equation}
In the case $\theta=2$, 
this generalization is suggested  by  the previous results for \eqref{eq:1.2}
in the framework of Sobolev spaces  $H^s_q$. 
Indeed, for any $s\in \R$ and $1<q<\infty$, 
the Cauchy problem \eqref{eq:1.2} has also been studied for initial data $\varphi\in H^s_q$, 
where
$$
H^s_q(\R^n):=\left\{\psi\in {\cal S}'(\R^n): (1-\Delta)^{\frac s2}\psi\in L^q (\R^n)\right \}
$$
(see, i.e., \cite{R}).
For the case $s<n/q$, 
a critical power nonlinearity appears in analogy with the theory in the Lebesgue framework. 
While for $s>n/q$,  no growth condition is necessary to establish the existence of a solution,
in the  case $s=n/q$ any power nonlinearity is allowed 
and  one wonders which is the  optimal critical growth at infinity.
By the Trudinger inequality in $\mathbb{R}^n$ we have the embedding
$$
H^s_q\hookrightarrow \exp L^{\Phi}.
$$
Here $\exp L^{\Phi}$ is the Orlicz space defined by the convex function
$$
\Phi(t):=\exp\left({t^{\frac q{q-1}}}\right )-\sum_{j=0}^{k_0-1}\frac {t^{jq/(q-1)}}{j!},
$$
where $k_0$ is the smallest integer satisfying $k_0\ge q-1$.
This indicates that,
for the Cauchy problem \eqref{eq:1.1} with $\varphi\in H^{n/q}_q$,
the critical growth of the nonlinearity at infinity should have the same rate 
as the case $\varphi\in\exp L^{q/(q-1)}$, that is, $f(u)\sim{e}^{u^{q/(q-1)}}$.

%
\section{Preliminaries}
In this section
we recall some properties of the fundamental solution of the fractional diffusion equation,
and give preliminary estimates.
Furthermore we prove the boundedness of the nonlinear term $f(u)$
under the smallness assumption for the solution $u$.

Let $\theta\in(0,2]$
and $G_\theta(t,x)$ be the fundamental solutions of the linear diffusion equations
\begin{equation}
\label{eq:3.1}
\partial_t u+{\cal L}_{\theta}u=0,
 \qquad t>0,\quad x\in{\mathbb R}^n,
\end{equation}
where ${\cal L}_{\theta}$ is given in \eqref{eq:2.12}.
It is well known that $G_\theta$ satisfies
the following (see, i.e., \cite{IKK1}):
\begin{itemize}
  \item[{\rm (i)}]
    $\displaystyle{\int_{{\mathbb R}^n}G_\theta(t,x)dx=1}$ for any $t>0$;\vspace{5pt}
  \item[{\rm (ii)}]
    For any $j\in\N_0$,
    $$
    |\nabla^jG_\theta(t,x)|\le Ct^{-\frac{n+j}{\theta}}\left(1+{t^{-1/\theta}|x|}\right)^{-(n+j)},
    \qquad t>0,\quad x\in{\mathbb R}^n.
    $$
    Furthermore, for any $1\le r\le\infty$,
    $$
    \sup_{t>0}\,t^{\frac{n}{\theta}(1-\frac{1}{r})}\|G_\theta(t)\|_{L^r}<\infty;
    $$
\item[{\rm (iii)}] For any $0<s\le t$,
    \begin{equation}
    \label{eq:3.2}
    G_\theta(t,x)=\int_{{\mathbb R}^n}G_\theta(t-s,x-y)G_\theta(s,y)\,dy.
    \end{equation}
\end{itemize}
For any $\varphi\in L^r$ ($1\le r\le\infty$),
we put
$$
{e}^{-t{\cal L}_\theta}\varphi(x):=\int_{{\mathbb R}^n}G_\theta(t,x-y)\varphi(y)\,dy, 
\qquad t>0,\quad x\in{\mathbb R}^n,
$$
which is a solution of \eqref{eq:3.1} with the initial data $\varphi$,
and it follows from \eqref{eq:3.2} that
$$
{e}^{-t{\cal L}_\theta}\varphi(x)={ e}^{-(t-s){\cal L}_\theta}[{ e}^{-s{\cal L}_\theta}\varphi](x),
\qquad t\ge s>0, \quad x\in{\mathbb R}^n.
$$
Combining property~(ii) with the Young inequality and \cite{IKK1,IKK2},
we have:
\begin{itemize}
\item[$(G_1)$] There exists a constant $c_\theta$, which depends only on $n$ and $\theta$,  such that
\begin{equation}
\label{eq:3.3} \|{ e}^{-t{\cal L}_\theta}\varphi\|_r\le
c_\theta t^{-\frac{n}{\theta}(\frac{1}{q}-\frac{1}{r})}\|\varphi\|_{L^q},
\qquad \|{ e}^{-t{\cal L}_\theta}\varphi\|_{L^q} \le \|\varphi\|_{L^q},
\qquad t>0,
\end{equation}
for any $\varphi\in L^q$ and $1\le q\le r\le\infty$;\vspace{3pt}
\item[$(G_2)$]
Let $\varphi\in L^1$ be such that
$$
\int_{{\mathbb R}^n}\varphi(x)\,dx=0.
$$
Then
$$
\lim_{t\to\infty}\|{e}^{-t{\cal L}_\theta}\varphi\|_{L^1}=0.
$$
\end{itemize}
We recall now the definition and the main properties of the Orlicz space ${\rm exp}L^r$ for $r>1$.
\begin{definition}
\label{Definition:3.1}
Let $r>1$.
We define the Orlicz space ${\rm exp}L^r$ as
\begin{equation*}
{\exp{L^r}}:=\left\{ u\in L^1_{loc}(\mathbb R^n), \exists \lambda>0 :
\int_{\mathbb R^n}\left({\exp}\left(\frac{|u(x)|}{\lambda}\right)^r-1\right)dx<\infty\right\},
\end{equation*}
\end{definition}
where the norm is given by the Luxemburg type
\begin{equation*}
\|u\|_{{\exp{L^r}}}:=\inf \left\{ \lambda>0 \ {\rm such \ that}\
\int_{\mathbb R^n}\left({\exp}\left(\frac{|u(x)|}{\lambda}\right)^r-1\right)dx\leq 1
\right\}.
\end{equation*}
The space ${\exp{L^r}}$ endowed with the  norm
$\|u\|_{{\exp{L^r}}}$ is a Banach space. Moreover
${\exp{L^r}}\hookrightarrow L^q$ for any $q\in[r,\infty)$.  
An example of a function that is not bounded and that belongs to ${\exp{L^r}}$ is
\begin{equation*}
u(x)=
\begin{cases}
(-\log|x|)^{1/r} &  |x|<1,\\
0 & |x|\geq 1.
\end{cases}
\end{equation*}
We stress also that  $C_0^\infty$  is not dense in ${\exp{L^r}}$. 
Finally the space ${\exp{L^r}}$ admits as predual the Orlicz space defined by the complementary function
of $A(t)=e^{t^r}-1$, denoted by $\tilde A(t)$. This  complementary function  is in particular
a convex function  such that 
$\tilde A(t)\sim t^2$ as $t\to 0$ and  $\tilde A(t)\sim t\log^{1/r}t$ as $t\to\infty$.
For more information about the Orlicz spaces, we refer to \cite{A,Rao}.
Furthermore the following estimates hold.
\begin{lemma}
\label{Lemma:3.1}
Let $r>1$.
Then, for any $p\in[1,r]$,
there exists a positive constant $C_\theta$, which depends only on $n$ and $\theta$, such that
\begin{align}
& 
\|{e}^{-t{\cal L}_\theta}\varphi\|_{{\exp}L^r}\le\|\varphi\|_{{\rm exp}L^r},\qquad t>0,
\notag
\\
&
\|{e}^{-t{\cal L}_\theta}\varphi\|_{{\exp}L^r} \le C_\theta t^{-\frac{n}{\theta p}}
\left(\log\left(t^{\frac{n}{\theta}}+1\right)\right)^{-\frac{1}{r}}\|\varphi\|_{L^p},
\qquad t>0.
\notag
\end{align}
\end{lemma}
A proof of  this lemma is based on the following basic estimates
and, for the case $r=2$, can be find it in \cite{I,I-t}.
\begin{lemma}
\label{Lemma:3.2}
Let $r>1$.
Then it holds that
\begin{equation}
\label{eq:3.4}
\|\psi\|_{L^p}\le\left[\Gamma\left(\frac{p}{r}+1\right)\right]^{\frac{1}{p}}\|\psi\|_{{\rm exp}L^r}
\end{equation}
for any $p\in[r,\infty)$,
where $\Gamma$ is the Gamma function
$$
\Gamma(q):=\int_0^\infty \xi^{q-1}e^{-\xi}\,d\xi, \qquad q>0.
$$
\end{lemma}
A proof of this lemma can be found in \cite{RT}.
\begin{lemma}
\label{Lemma:3.3}
For any $p\ge1$ and $r\ge1$,
there exists a positive constant $C$, which is independent of $p$ and $r$,
such that
\begin{equation}
\label{eq:3.5}
\Gamma(rp+1)^{\frac{1}{p}}\le C\Gamma(r+1)p^{r}.
\end{equation}
\end{lemma}
\noindent {\bf Proof.}
This is a consequence of Stirling's formula. 
Indeed, for  a fixed  positive constant  $C>1$, 
there exists $r_0$ large enough such that, 
for any $r\geq r_0$ and $p\geq 1$,
$$
\frac{{\Gamma(rp+1)}^{1/p}}{{\Gamma(r+1)}}
\leq C\frac{(rp)^{r}{\rm e}^{-r}(2\pi rp)^{1/(2p)}}{r^{r}{\rm e}^{-r}(2\pi r)^{1/2}}\leq Cp^{r}.
$$
For $1\leq r\leq r_0$, 
we consider first the case of large values of $p$, namely $p\geq p_0$,  
and we obtain an estimate similar to the previous one.
Finally we observe that for  $1\leq r\leq r_0$ and  $1\leq p\leq p_0$ the quotient is bounded by a constant.
$\Box$\vspace{7pt}

At the end of this section we give the following estimate on the nonlinearity $f(u)$,
which is crucial throughout this paper.

\begin{lemma}
\label{Lemma:3.4}
Let $n\geq 1$ and $M>0$.
Suppose that the function $u\in L^\infty(0,\infty\,;\, \exp L^2)$ satisfies the condition
$$
\sup _{t>0} \|u(t)\|_{\exp{L^2}} \leq M.
$$
Let $f$ be the function defined as in \eqref{eq:1.6}.
Furthermore, let $p_1$ be the constant given in \eqref{eq:2.7}, namely $p_1= \max \left(1, (2n)/(4+n)\right )$.
Then, for all $p\in [p_1,\infty)$, there exists $\eps =\eps (p) >0$ such that, 
if $M<\eps$, then
\begin{equation}
\label{eq:3.6}
\sup_{t>0} \|f(u(t))\|_{L^p} \leq 2Cp^3M^{1+\frac{4}{n}},
\end{equation}
where $C$ is independent of $p$, $n$ and $M$.
\end{lemma}
\noindent{\bf Proof.}
For any $k\in{\mathbb N}_0$, we put
\begin{equation}
\label{eq:3.7}
\ell_k:=2k+1+\frac{4}{n}.
\end{equation}
Then,
since it holds from $p\ge p_1$ with \eqref{eq:2.7} that
$$
\ell_kp\ge \left(1+\frac{4}{n}\right)\frac{2n}{n+4}=2
$$
for any $k\in{\mathbb N}_0$,
by \eqref{eq:1.6} and \eqref{eq:3.4}
we have
\begin{equation}
\label{eq:3.8}
\begin{split}
\|f(u(t))\|_{L^p}
&
\leq \sum_{k=0}^{\infty}
\frac {\|u^{\ell_k}(t)\|_{L^p}}
{k!}
\\
&
 =
 \sum_{k=0}^{\infty} \frac {\|u(t)\|^{\ell_k}_{L^{\ell_kp}}}
{k!}
\leq
\sum_{k=0}^{\infty}
\frac {\left( \Gamma (\frac{\ell_kp}{2}+1)^{\frac{1}{\ell_kp}} \|u(t)\|_{\exp L^2} \right )^{\ell_k}}{k!}
\end{split}
\end{equation}
for all $t>0$.
By \eqref{eq:3.5} with the monotonicity property of the Gamma function we see that
\begin{equation*}
\begin{split}
 \Gamma \left(\frac{\ell_kp}{2}+1\right)^{\frac{1}{p}}
&
 \le C\Gamma\left(\frac{\ell_k}{2}+1\right)p^{\frac{\ell_k}{2}}
 \\
 &
 = C\Gamma\left(k+\frac{3}{2}+\frac{2}{n}\right)p^{\frac{\ell_k}{2}}
 \\
 &
 \le C\Gamma(k+4)p^{\frac{\ell_k}{2}}
 = C (k+3)! p^{\frac{\ell_k}{2}}.
 \end{split}
\end{equation*}
This together with the assumption on $u$ and \eqref{eq:3.8} implies that
\begin{equation*}
\begin{split}
\|f(u(t))\|_{L^p}
&
\le    C \sum_{k=0}^{\infty} \frac{(k+3)!}{k!}\left(p\|u(t)\|_{\exp L^2}^2\right)^{\frac{\ell_k}{2}}
\\
&
\le  C \sum_{k=0}^{\infty}(k+1)(k+2)(k+3)\left(pM^2\right )^{k+\frac{1}{2}+\frac{2}{n}},
\qquad t>0.
\end{split}
\end{equation*}
Now for $M <\eps (p)$ small enough we get
$$
 \sup_{t>0} \|f(u(t))\|_{L^p} \leq C\frac {(pM^2)^{\frac{1}{2}+\frac{2}{n}}}{(1-pM^2)^4}
 \leq 2C(pM^2)^{\frac{1}{2}+\frac{2}{n}} \leq 2Cp^3M^{1+\frac{4}{n}}.
$$
This implies \eqref{eq:3.6},
and the proof of Lemma~\ref{Lemma:3.4} is complete.
$\Box$
\vspace{7pt}

%
\section{Equivalence and uniqueness}
In this section we prove 
 the equivalence between small weak and 
small mild solutions of \eqref{eq:1.1}.
Furthermore, we show that the small solution is unique in any dimension $n\in \N$.
\vspace{5pt}

We first prove Proposition~\ref{Theorem:2.1}.

\noindent{\bf Proof of Proposition~\ref{Theorem:2.1}.}
Let us first remark that, for $\varepsilon $ small enough, 
by \eqref{eq:1.1} we can apply Lemma \ref{Lemma:3.4},
and we see that $f(u)\in L^{\infty}(0,\infty; L^{p_1})$ where $p_1$ is the constant given in \eqref{eq:2.7}.

Suppose that  $u$ is a weak solution of the Cauchy problem \eqref{eq:1.1}.
Then, for $0<s<t$, it holds 
\[
{e}^{(t-s)\Delta}f(u(s))
= {e}^{(t-s)\Delta}\left(\partial_s u-\Delta u\right)= \partial_s\left({e}^{(t-s)\Delta}u(s)\right)
\quad\mbox{in}\quad{\mathcal D}'.
\]
Integrating on $(\tau,t)$, we obtain
$$
u(t)={e}^{(t-\tau)\Delta}u(\tau)+\int_{\tau}^t{e}^{(t-s)\Delta}f(u(s))\,ds
$$
for all $t>\tau>0$.
Let $\tau\to 0$.
Since $f(u)\in L^{\infty}(0,T; L^{p_1})$, 
we get
$$
\int_{\tau}^t{e}^{(t-s)\Delta}f(u(s))\,ds\to \int_{0}^t{e}^{(t-s)\Delta}f(u(s))\,ds
\quad\mbox{in}\quad L^{p_1}.
$$
Moreover, it follows from $u(\tau)\buildrel{*}\over\rightharpoonup\varphi$ in $\exp L^2$ that
$$
{e}^{(t-\tau)\Delta}u(\tau)\buildrel{*}\over\rightharpoonup{e}^{t\Delta}\varphi,
$$
and the same limit holds in ${\cal D}'$.
So we obtain
$$
u(t)={e}^{t\Delta}\varphi+\int_{0}^t{e}^{(t-s)\Delta}f(u(s))\,ds
\quad\mbox{in}\quad\mathcal{D}'.
$$
Moreover ${e}^{t\Delta}\varphi \in C((0,\infty); \exp L^2)$ and by Lemma \ref{Lemma:3.4} the integral term
belongs to $C([0,\infty), \exp L^2)$. This completes the proof that 
statement $($i$)$ implies statement $($ii$)$.

Suppose now that $u\in L^{\infty}(0,T; \exp L^2)$ is a mild solution of \eqref{eq:1.1}.
We know that
$$
\partial_t{e}^{t\Delta}\varphi=\Delta{e}^{t\Delta }\varphi \quad \text{in}\quad{\mathcal D}'.
$$
Let us consider $\eta\in {\cal D}((0,T)\times \R^n)$. 
Of course, there exists $\tilde T\in(0,T)$ such that $\eta (t,x)=0$ for $t\in[\tilde T,T)$.
By a direct computation we are going to prove that
\begin{equation*}
\begin{split}
&
\langle\partial_t\int_0^te^{(t-s)\Delta}f(u(s))\,ds,\eta(t)\rangle_{L^2(L^2)}
\\
&
=\langle f(u),\eta\rangle_{L^2(L^2)}+\langle\Delta\int_0^t{e}^{(t-s)\Delta}f(u(s))\,ds,\eta(t)\rangle_{L^2(L^2)}.
\end{split}
\end{equation*}
Indeed,
\begin{equation*}
\begin{split}
&\langle \partial_t\int_0^t{e}^{(t-s)\Delta}f(u(s))\,ds ,\eta(t)\rangle_{L^2(L^2)}\\
&=-\langle\int_0^t{e}^{(t-s)\Delta}f(u(s))\,ds ,\partial_t\eta(t)\rangle_{L^2(L^2)}\\
&=-\int_0^{\tilde T}\langle\int_0^t{e}^{(t-s)\Delta}f(u(s))\,ds ,\partial_t\eta(t)\rangle_{L^2}\,dt\\
&=-\int_0^{\tilde T}\int_0^t\langle{e}^{(t-s)\Delta}f(u(s)),\partial_t\eta(t)\rangle_{L^2}\,ds\,dt\\
&=-\int_0^{\tilde T}\int_0^t\langle f(u(s)),{e}^{(t-s)\Delta}\partial_t\eta(t)\rangle_{L^2}\,ds\,dt\\
&=\int_0^{\tilde T}\int_0^t\langle f(u(s)),-\partial_t({e}^{(t-s)\Delta}\eta(t))
+\Delta {e}^{(t-s)\Delta}\eta(t)\rangle_{L^2}\,ds\,dt\\
&=\int_0^{\tilde T}\int_0^t-\partial_t\langle f(u(s)),{e}^{(t-s)\Delta}\eta(t)\rangle_{L^2}\,ds\,dt
+\int_0^{\tilde T}\int_0^t\langle f(u(s)),\Delta {e}^{(t-s)\Delta}\eta(t)\rangle_{L^2}\,ds\,dt\\
&=\int_0^{\tilde T}\int_s^{\tilde T}-\partial_t\langle f(u(s)),{e}^{(t-s)\Delta}\eta(t)\rangle_{L^2}\,dt\,ds
+\int_0^{\tilde T}\langle\int_0^t\Delta {e}^{(t-s)\Delta}f(u(s))\,ds ,\eta(t)\rangle_{L^2}\,dt\\
&=\int_0^{\tilde T}\langle f(u(s)),\eta(s) \rangle_{L^2}\,ds 
+\langle\Delta\int_0^t {e}^{(t-s)\Delta}f(u(s))\,ds ,\eta(t)\rangle_{L^2(L^2)} \\
&=\langle f(u),\eta\rangle_{L^2(L^2)}
+\langle\Delta\int_0^t {e}^{(t-s)\Delta}f(u(s))\,ds ,\eta(t)\rangle_{L^2(L^2)}.
\end{split}
\end{equation*}
Therefore we  obtain
$$
\partial_tu=\Delta u+f(u) \quad{\rm in}\quad{\mathcal   D}'.
$$
This yields that statement $($ii$)$ implies statement $($i$)$,
and the proof of Theorem~\ref{Theorem:2.1} is complete.
$\Box$\vspace{7pt}

Next we show the uniqueness of the small solution of \eqref{eq:1.1}.
\medskip

\noindent{\bf Proof of Proposition~\ref{uniqueness}.}
Let us suppose that $u,v$ are two solutions of the Cauchy problem \eqref{eq:1.1}.
Furthermore, since they converge to $\varphi$ weak* in $\exp L^2$,
we can continue them in $t=0$ as $\varphi$.
Moreover, by Lemma \ref{Lemma:3.4} we have
$$
\|u(t)-v(t)\|_{\exp L^2}
\le
\|u(t)-e^{t\Delta}\varphi\|_{\exp L^2}
+
\|v(t)-e^{t\Delta}\varphi\|_{\exp L^2}\to0
\quad\mbox{as}\quad t\to+0.
$$
as $t\to+0$.
Let
\be
\label{eq:4.1}
H(t):=\|u(t)-v(t)\|_{\exp L^2}, \quad t\geq 0,
\ee
so that $H(0)=0$.
Furthermore, put
$$
t_0:=\sup \left\{ t\in [0,\min(T,T'))\ \text { such that }
H(s)=0\   \text{ for  every } s\in [0,t]\right \}.
$$
By \eqref{eq:4.1} we have $t_0\geq 0$.
 By
contradiction we assume  $t_0<\min (T,T')$. Since $H(t)$
is continuous in time  we have $H(t_0)=0$. Let us denote
$\tilde u(t)=u(t+t_0)$ and $\tilde v(t)=v(t+t_0)$ so that $\tilde u$
and $\tilde v$ satisfy equation \eqref{eq:1.1} on
$(0,\infty)$ and by $H(t_0)=0$ it follows $\tilde u(0)=\tilde v(0)$. We will prove
that there exists a positive time $\tilde t $ such that
\begin{equation}
\label{eq:4.2}
\sup_{0<t<\tilde t}\|\tilde u(t)-\tilde v(t)\|_{{\exp}L^2}\leq
C(\tilde  t)\sup_{0<t<\tilde t}\|\tilde u(t)-\tilde v(t)\|_{{\exp}L^2}
\end{equation}
for a constant $C(\tilde t)<1$, and so $\| \tilde u(t)-\tilde v(t)\|_{\exp L^2}=0$
for any $t\in [0,\tilde t]$. 
Therefore $u(t+t_0)=v(t+t_0)$ for any $t\in [0,\tilde t] $ in contradiction with the definition of $t_0$.
In order to establish inequality \eqref{eq:4.2} we control both
the $L^2$ norm and  the $L^\infty$ norm of the difference of the two solutions.
By \eqref{eq:1.6} we see that there exist positive constants $ C$ and $ \lambda$ such that
\[
\left|f(x)-f(y) \right|\leq  C |x-y| \left({\rm e}^{ \lambda x^2}+{\rm e}^{ \lambda y^2}\right),
\qquad x,y\in\mathbb R.
\]
Then, for $p$, $q>2$ such that $1/p+1/q=1/2$, we have
\[
\begin{split}
\|\tilde u(t)-\tilde v(t)\|_{L^2} 
&\leq  C\int_0^t \left\| |\tilde u(s)-\tilde v(s)| 
({e}^{ \lambda \tilde u^2(s)} +{ e}^{ \lambda \tilde v^2(s)}) \right\|_{L^2}\,ds
\\
&\leq C\int_0^t \| \tilde u(s)-\tilde v(s)\|_{L^2}\,ds
\\
& \qquad +  C\int_0^t \left(\| \tilde u(s)-\tilde v(s)\|_{L^q}
\left\|({e}^{ \lambda \tilde u^2(s)}-1)+({e}^{ \lambda \tilde v^2(s)}-1) \right\|_{L^p} \right)\,ds
\\
&\leq Ct\sup_{0<s<t}\|\tilde u(s)-\tilde v(s)\|_{\exp L^2}
\\
& \qquad + C \sup_{0<s<t}\|\tilde u(s)-\tilde v(s)\|_{\exp L^2}
\int_0^t \left\|({ e}^{ \lambda \tilde u^2(s)}-1) +(e^{ \lambda \tilde v^2(s)}-1) \right\|_{L^p}\,ds
\end{split}
\]
for all $t>0$.
Moreover, 
if $\sup_{t>0}\|\tilde u(t)\|_{\exp L^2}$ and $\sup_{t>0}\|\tilde v(t)\|_{\exp L^2}$ are small enough,
then, by the same proof as Lemma~\ref{Lemma:3.4} 
the term in the integral is uniformly bounded in time.
Indeed,
\begin{equation}
\label{eq:4.3}
\begin{split}
\sup_{0<s<\infty} & \left\|({e}^{\tilde \lambda \tilde u^2(s)}-1) +
({e}^{\tilde \lambda \tilde v^2(s)}-1) \right\|_{L^p}
\\
&\leq C(\tilde u,\tilde v)<\infty.\\
\end{split}
\end{equation}
Therefore, we obtain
\begin{equation}
\label{eq:4.4}
\sup_{0<s<t}\|\tilde u(s)-\tilde v(s)\|_{L^2} \leq C(\tilde u,\tilde v)\,
t \sup_{0<s<t}\|\tilde u(s)-\tilde v(s)\|_{{\exp}L^2}.
\end{equation}
In a similar way, for $r> 2$ such that $n/(2r) <1$, we get
\[
\begin{split}
\|\tilde u(t)-\tilde v(t)\|_{L^\infty} 
&\leq C\int_0^t(t-s)^{-\frac{n}{2r}} \left\| |\tilde u(s)-\tilde v(s)|
  \left({e}^{ \lambda \tilde u^2(s)}+{e}^{ \lambda \tilde v^2(s)}
  \right)
\right\|_{L^r}\,ds
\\
&\leq C\int_0^t (t-s)^{-\frac n{2r}}\|\tilde u(s)-\tilde v(s)\|_{L^r}\,ds
\\
&\qquad + C \int_0^t(t-s)^{-\frac{n}{2r}}\|\tilde u(s)-\tilde v(s)\|_{L^{\bar q}} 
\left\|{e}^{ \lambda \tilde u^2(s)}-1+
{e}^{ \lambda \tilde v^2(s)}-1
\right\|_{L^{\bar p}}\,ds
\end{split}
\]
for some $\bar q$, $\bar p$ such that $1/{\bar q}+1/{\bar p}=1/r$. 
Since $ \bar p \ge r >2$, 
we can apply  an estimate similar to \eqref{eq:4.3},  
and obtain that
\begin{equation}
\label{eq:4.5}
\begin{split}
&\sup_{0<s<t}\|\tilde u(s)-\tilde v(s)\|_{L^{\infty}} 
\leq C(\tilde u,\tilde v)\,t^{1-\frac{n}{2r}}\,\sup_{0<s<t}\|\tilde u(s)-\tilde v(s)\|_{{\exp}L^2}.\\
\end{split}
\end{equation}
Therefore the two inequalities \eqref{eq:4.4} and \eqref{eq:4.5} imply
$$
\sup_{0<s<t}\|\tilde u(s)-\tilde v(s)\|_{{\exp}L^2}
\leq C(\tilde u,\tilde v) (t^{1-\frac{n}{2r}}+t)
\sup_{0<s<t}\|\tilde u(s)-\tilde v(s)\|_{{\exp}L^2} \, ,
$$
and, for $t$ small enough, 
we obtain the desired estimate.  
$\Box$\vspace{7pt}

\section{Singular initial data}

In order to prove Theorem \ref{Theorem:2.2} the principal idea is the following.
Even if the initial data $\varphi$ satisfy  $\varphi \in \exp L^2$
and so  $\varphi \in L^p$ for $p\in [2,\infty)$ with a norm blowing up with $p$,
the solution is indeed in $L^2\cap L^\infty$ for all $t>0$.
Moreover, for each fixed $t_0>0$, 
the norm  $\|u(t_0)\|_\infty$ is arbitrarily small provided the initial data are sufficiently small in $\exp L^2$.
So, thanks to the uniqueness result (Proposition~\ref{uniqueness}),
it is possible to consider the solution $u(t)$ for $t\geq t_0$
as the limit of a recursive procedure building up a solution starting
with initial data belonging to all $L^p$ with $p\in [2,\infty]$
and with norm uniformly bounded.
\vspace{5pt}

Let us prove now \eqref{eq:2.2} for small time.
\begin{lemma}\label{Lemma:5.1}
Let $n\geq 1$.
There is $\eps = \eps (n) >0$ such that, 
if $\|\varphi\|_{\exp L^2} <\eps$, 
then the unique solution $u$   of  the Cauchy problem \eqref{eq:1.1} satisfying \eqref{eq:1.8}
and
\be
\label{eq:5.1}
\|u(t)\|_{L^q} \leq C t^{-\frac n2 \left( \frac 12 -\frac 1q\right )}\|\varphi\|_{\exp L^2},\quad 0<t\leq 1,
\ee
for any $q\in[2,\infty]$,
where $C>0$ depending only on $n$.
\end{lemma}
\noindent{\bf Proof.}
Let us consider
\[
u(t)= e^{t\Delta} \varphi +  \int_0^t e^{(t-s)\Delta} f(u(s)) \, ds := e^{t\Delta} \varphi + F(u(t))
\]
and let us denote as in Lemma \ref{Lemma:3.4}
\[
M = 2\|\varphi\|_{\exp L^2}.
\]
As for the linear part,
for any $q\in[2,\infty]$,
by \eqref{eq:3.3} and \eqref{eq:3.4} we have
\[
\|e^{t\Delta} \varphi \|_{L^q} 
\leq c_2 t^{-\frac n2\left(\frac 12 -\frac 1q\right )}\|\varphi \|_{L^2}
\leq  c_2 t^{-\frac n2\left(\frac 12 -\frac 1q\right )}\|\varphi \|_{\exp L^2},\quad t>0.
\]
Let us consider now the the nonlinear part $F(u)$.
By choosing $r\geq p_1$ such that $2r>n$ (and so $r$ depends only on the dimension $n$),
we get
\[
\|F(u(t))\|_{L^\infty} 
\leq c_2 \int_0^{t} (t-s)^{-\frac n{2r}} \|f(u(s))\|_{L^r} \, ds
\leq C t^{1-\frac n{2r}} \sup_{t>0}\|f(u(t))\|_{L^r},
\qquad t>0.
\]
Now, by Lemma \ref{Lemma:3.4} we get for $M$ small enough depending on $r$ (and so on $n$ only)
\[
C t^{1-\frac n{2r}} \sup_{t>0}\|f(u(t))\|_{L^r} \leq  C r^3 t^{1-\frac n{2r}} M^{1+\frac{4}{n}}.
\]
Now, for $t\in (0,1]$ and $M$ small enough, depending again only on $n$, we get
\[
\sup_{0<t\leq 1} \|F(u(t))\|_{L^\infty} \leq {C} \|\varphi\|_{\exp L^2}.
\]
So we get 
\be
\label{eq:5.2}
\|u(t)\|_{L^\infty} 
\leq c_2t^{-\frac n4} \|\varphi\|_{\exp L^2} + C \|\varphi\|_{\exp L^2}  
\leq  C  t^{-\frac n4} \|\varphi\|_{\exp L^2}, \qquad t\in (0,1].
\ee
On the other hand, by the embedding $\exp L^2 \subset L^2$ \eqref{eq:3.4} we get
$$
\|u(t)\|_{L^2} \leq \Gamma (2)^\frac 12 \|u(t)\|_{\exp L^2}, \qquad   t>0.
$$
This together with \eqref{eq:5.2} implies \eqref{eq:5.1}, 
and the proof of Lemma~\ref{Lemma:5.1} is complete.
$\Box$\vspace{7pt}

%
We have now to prove \eqref{eq:2.2} for large times $t>1$.
Let us put
\[
v(t,x):= u(t+1,x)
\]
so that
\[
\left\{
\begin{aligned}
&\partial_t v= \Delta v + f(v), \quad t>0,\  x\in\R^n \\
&v(0,x)=u(1,x) \in L^q, \quad x\in \R^n, \quad  q\in [2,\infty],
\end{aligned}
\right .
\]
and by Lemma \ref{Lemma:5.1}  for all $q\in [2,\infty]$
we have $\|u(1)\|_{L^q} \leq  C \|\varphi\|_{\exp L^2} \leq \eps$ small enough.
Let us consider now the following recursive sequence $\{v_j\}_{j\in \N} $
\be\label{eq:5.3}
\left\{
\begin{aligned}
&v_0(t)= e^{t\Delta} v(0)\\
&v_{j+1}(t)= v_0(t) +  \int_0^t e^{(t-s)\Delta} f(v_j(s))\,ds, \quad j\geq 0.
\end{aligned}
\right .
\ee
Since $0\leq v_j(t,x) \leq v_{j+1} (t,x) \leq v(t,x)$ for $x\in\R^n$ and $t> 0$,  we get by uniqueness
\be\label{eq:5.4}
\lim_{j\to \infty} v_j(t,x)= v(t,x) =u(t+1,x)
\ee
the limit being considered pointwise.
We   prove now \eqref{eq:2.2} for large times through the recursive sequence \eqref{eq:5.3}.
We begin by a Lemma  which describes how boundedness and decreasing in time propagate on the nonlinear term.

\begin{lemma}
\label{Lemma:5.2}
Let $n\geq 1$. Let $M=2\|\varphi\|_{\exp L^2}$ and  $v(x,t)$ a function satisfying
\begin{equation}
\label{eq:5.5}
\sup_{s>0} (1+s)^{\frac{n}{2}\left(\frac{1}{2}-\frac{1}{q}\right )}\|v(s)\|_{L^q} \leq C M,\qquad q\in[2,\infty],
 \end{equation}
 where $C$  independent of $q$.
Then,  there exists $\eps >0$ such that, 
if $M<\eps$, then, for any $r\in [p_1,\infty]$,
\begin{equation}
\label{eq:5.6}
\sup_{s>0}  (1+s)^{\frac{n}{2}\left(\frac{1}{2}-\frac{1}{r}\right )+1} \|f(v(s))\|_{L^r} \leq 2 (CM)^{1+\frac{4}{n}}\end{equation}
where $f$ and $p_1$ are defined as in \eqref{eq:1.6} and \eqref{eq:2.7}, respectively.
\end{lemma}
\noindent{\bf Proof.}
Let $k\in{\mathbb N}_0$ and $\ell_k$ be the constant given in \eqref{eq:3.7}.
Then, for any $r\in[p_1,\infty]$, by \eqref{eq:1.6} and \eqref{eq:5.5} we have
\begin{equation}
\label{eq:5.7}
\begin{split}
\|f(v(s))\|_{L^r}
&
\le \sum_{k=0}^{\infty} \frac 1{k!} \|v(s)\|_{L^{\ell_kr}}^{\ell_k}
\\
&
\le  \sum_{k=0}^{\infty} \frac 1{k!} \left( C (1+s)^{-\frac{n}{2} \left( \frac{1}{2}
-\frac{1}{\ell_kr}\right )} M\right )^{\ell_k}
\\
&
\le (CM)^{1+\frac{4}{n}}(1+s)^{\frac{n}{2r}-\frac{n}{4}\left(1+\frac{4}{n}\right)}
 \sum_{k=0}^{\infty} \frac 1{k!} \left( C (1+s)^{-\frac{n}{4}}M\right)^{2k}
 \\
 &
 \le (CM)^{1+\frac{4}{n}}(1+s)^{-\frac{n}{2}\left(\frac{1}{2}-\frac{1}{r}\right)-1}
 \sum_{k=0}^{\infty} \frac 1{k!} (CM)^{2k}
\end{split}
\end{equation}
for all $s>0$.
We can take a sufficiently small $\varepsilon$, which is  independent of $r$, so that,
for $M\leq \eps $,
it holds that
$$
\sum_{k=0}^{\infty} \frac 1{k!} \left( CM\right)^{2k}=e^{(CM)^2}\le 2.
$$
This together with \eqref{eq:5.7} implies \eqref{eq:5.6}.
Thus Lemma~\ref{Lemma:5.2} follows.
$\Box$\vspace{7pt}

Now we are in position to prove Theorem~\ref{Theorem:2.2}.
\vspace{5pt}

\noindent{\bf Proof of Theorem \ref{Theorem:2.2}.}
By Lemma~\ref{Lemma:5.1} we have \eqref{eq:2.2} for $0<t\le 1$.
So it suffices to prove \eqref{eq:2.2} for $t>1$.
Let $j\in {\mathbb N}_0$ and $v_j(t,x)$ be the function given in \eqref{eq:5.3}.
Since  $v_0(t) =e^{t\Delta} v(0) = e^{t\Delta} u(1)$,
for any $q\in[2,\infty]$,
we have by \eqref{eq:3.3}
\[
\|v_0(t)\|_{L^q} \le c_2 t^{-\frac{n}{2}\left( \frac{1}{2} -\frac{1}{q}\right)}\|v(0)\|_{L^2}, \quad t>0,
\]
and also
\[
\|v_0(t)\|_{L^q} \le \|v(0)\|_{L^q}, \quad t>0.
\]
On the other hand, by Lemma~\ref{Lemma:5.1} we have
\[
\| v(0)\|_{L^q}= \|u(1)\|_{L^q} \le C \|\varphi\|_{\exp L^2}.
\]
Therefore, for any $q\in [2,\infty]$,
there exists a positive constant $c_*$ depending only on $n$ such that
\begin{equation}
\label{eq:5.8}
\|v_0(t) \|_{L^q} \le c_* (1+t)^{-\frac{n}{2}\left( \frac{1}{2} -\frac{1}{q}\right)}  \|\varphi\|_{\exp L^2}.
\end{equation}
Let $m\in{\mathbb N}_0$, and put
$$
\tilde{M}:=c_*\|\varphi\|_{\exp L^2}.
$$
Let us assume that, for $j=m$,
\begin{equation}
\label{eq:5.9}
\|v_m(t)\|_{L^q} \leq 2\tilde M (1+t)^{-\frac{n}{2} \left( \frac{1}{2} -\frac{1}{q}\right )},\quad t>0,
\end{equation}
and let us prove it for $j=m+1$.
By \eqref{eq:5.3} it suffices to consider the nonlinear term.
For any $q\in [2,\infty]$, we put
\begin{equation*}
\begin{split}
&
\left\|  \int_0^t e^{(t-s)\Delta} f(v_m(s))\,ds \right \|_{L^q}
\\
 &\le
  \int_0^{t/2}\left\| e^{(t-s)\Delta} f(v_m(s))\right \|_{L^q}\, ds
  + \int_{t/2}^{t}\left\| e^{(t-s)\Delta} f(v_m(s))\right \|_{L^q}\, ds\\
 & =: J_1(t) +J_2(t),
 \qquad t>0.
  \end{split}
\end{equation*}
For the term $J_1(t)$,
by \eqref{eq:3.3} and \eqref{eq:5.6} we obtain
\begin{equation}
\label{eq:5.10}
\begin{split}
J_1(t)
&
\le   \int_0^{t/2}\left\| f(v_m(s))\right \|_{L^q}\, ds
\\
&
\le 2C{\tilde M}^{1+\frac{4}{n}}\int_0^1(1+s)^{-\frac{n}{2}(\frac{1}{2}-\frac{1}{q})-1}\, ds
\le  2 C {\tilde M}^{1+\frac{4}{n}}
\end{split}
\end{equation}
for all $0<t\le2$.
Moreover, since $1\le p_1<2$, by \eqref{eq:3.3} and \eqref{eq:5.6} we have
\begin{equation*}
\begin{split}
J_1(t)
&\le  c_2  \int_0^{t/2} (t-s)^{-\frac{n}{2}\left(\frac{1}{p_1}-\frac{1}{q}\right )}\|f(v_m(s))\|_{L^{p_1}} ds
\\
&
\le C{\tilde M}^{1+\frac{4}{n}}  t^{-\frac{n}{2}\left(\frac{1}{p_1}-\frac{1}{q}\right )}
\int_0^{\frac{t}{2}}(1+s)^{-\frac{n}{2}(\frac{1}{2}-\frac{1}{p_1})-1}ds
\\
&
\le C{\tilde M}^{1+\frac{4}{n}}
t^{-\frac{n}{2}\left(\frac{1}{p_1}-\frac{1}{q}\right )}(1+t)^{\frac{n}{2}(\frac{1}{p_1}-\frac{1}{2})}
\le
C{\tilde M}^{1+\frac{4}{n}}  t^{-\frac{n}{2}\left(\frac{1}{2}-\frac{1}{q}\right )}
\end{split}
\end{equation*}
for all $t\ge1$.
This together with \eqref{eq:5.10} implies that
\begin{equation*}
J_1(t)
\le D_1 {\tilde M}^{1+\frac{4}{n}} (1+t)^{-\frac{n}{2}\left(\frac{1}{2} -\frac{1}{q}\right )},
\qquad t>0,
\end{equation*}
where $D_1$ is a positive constant  independent of $m$ and $t$.
On the other hand,
for the term $J_2(t)$,
exploiting \eqref{eq:3.3} and \eqref{eq:5.6} again, we see that
\begin{equation*}
\begin{split}
J_2(t)
\le \int_{t/2}^{t}\left\| f(v_m(s))\right \|_{L^q}\, ds
&
\le 2C{\tilde M}^{1+\frac{4}{n}}
\int_{t/2}^t(1+s)^{-\frac{n}{2}\left(\frac{1}{2}-\frac{1}{q}\right )-1}\, ds
\\
&
\le  D_2{\tilde M}^{1+\frac{4}{n}}(1+t)^{-\frac{n}{2}\left(\frac{1}{2}-\frac{1}{q}\right )},
\qquad t>0,
\end{split}
\end{equation*}
where $D_2$ is a positive constant which independent of $m$ and $t$.
We can now assume
$$
(D_1+D_2){\tilde M}^{\frac{4}{n}}=(D_1+D_2)(c_*\|\varphi\|_{\exp L^2})^{\frac{4}{n}}<1
$$
in order to get
\begin{equation}
\label{eq:5.11}
\left\|  \int_0^t e^{(t-s)\Delta} f(v_m(s))\,ds \right \|_{L^q}
\le J_1(t)+J_2(t)
\le \tilde M (1+t)^{\frac n2\left(\frac 12-\frac 1q\right )},\qquad t>0.
\end{equation}
By \eqref{eq:5.3}, \eqref{eq:5.8} and \eqref{eq:5.11} we get the estimate \eqref{eq:5.9} on $v_{m+1}$.
Therefore we have the estimate \eqref{eq:5.9} for all $m\in{\mathbb N}_0$.
Passing to the limit with \eqref{eq:5.4} we see that
$$
\|u(t+1)\|_{L^q}
=
\|v(t)\|_{L^q} \le 2c_*(1+t)^{-\frac{n}{2} \left( \frac{1}{2} -\frac{1}{q}\right )}\|\varphi\|_{\exp L^2},\quad t>0.
$$
This implies \eqref{eq:2.2} for $t>1$,
and the proof of Theorem~\ref{Theorem:2.2} is complete.
$\Box$\vspace{7pt}
%
\section{Regular initial data}

We are going to prove now Theorem \ref{Theorem:2.3}.
For the estimate \eqref{eq:2.5} for small times,
we are going to exploit again Lemma \ref{Lemma:3.4} but this time,
in order not to assume smallness on  $\|\varphi\|_{L^p}$,
we introduce a time $T_*$, to be suitably chosen,
and we will split the estimate $\eqref{eq:2.5}$ into $0<t\leq 2T_*$ and $t>2T_*$.
The same splitting would have worked also in the singular case
without leading however to any significative improvement.
Let
$$
M:=\|\varphi\|_{\exp L^2}.
$$
By Proposition~\ref{Proposition:1.1}, it is not restrictive to assume $M<1$.
For  $\|\varphi\|_{p_*}>0$  we have
\[
M \le 2\|\varphi\|_{\exp L^2\cap L^{p_*}} \le 2\max\{2M,2\|\varphi\|_{L^{p_*}}\} 
\le 4\max\{1,\|\varphi\|_{L^{p_*}}\}.
\]
Let now
\be\label{eq:6.1}
K:= 4\max\{1,\|\varphi\|_{L^{p_*}}\}
\ee
and so
\be\label{eq:6.2}
M\leq K.
\ee
Then we begin by a Lemma analogous to Lemma \ref{Lemma:5.2}.
\begin{lemma}
\label{Lemma:6.1}
Let $n\geq 1$ and $p\in [1,2)$.
Furthermore let $p_*$ be the constant given in \eqref{eq:2.3}.
Suppose that $v(t,x)$ is a function satisfying
$$
\sup_{s>0}s^{\frac{n}{2}(\frac{1}{p_*} -\frac{1}{q})}\|v(s)\|_{L^q} \leq C K
$$
for any $q\in[p_*,\infty]$,
where $C$ is independent of $q$.
Then, there exists a sufficiently large constant $T_1= T_1(K,p_*)>1$ such that,
for any $r\in [p_2,\infty]$,
\begin{equation}
\label{eq:6.3}
\sup_{s>T_1}  s^{\frac{n}{2}(\frac{1}{p_*}-\frac{1}{r})+\frac{2}{p_*}}
\|f(v(s))\|_{L^r} \leq 2 (CK)^{1+\frac{4}{n}},
\end{equation}
where $f$ is defined as in \eqref{eq:1.6} and
\begin{equation}
\label{eq:6.4}
p_2=\max\left\{\frac{p_*n}{n+4},1\right\}.
\end{equation}
\end{lemma}
\begin{remark}
\label{Remark:6.1}
Comparing with Lemma~$\ref{Lemma:5.2}$,
since we are not assuming boundedness of $\|v(s)\|_{L^q}$ near $s=0$,
we can only obtain boundedness and decreasing of the nonlinear term for large times.
\end{remark}
\begin{remark}
For $1\leq n\leq 4$, since $p\in [1,2)$ and $(2n)/(n+4) \leq 1$, we have $p_*=p$ and $p_2=1$.
On the contrary for $n\geq 5$,  then $p_*$ might be  strictly greater than $p$ and $p_2$ might be strictly greater than $1$.
\end{remark}
\noindent{\bf Proof.}
Let $\ell_k=2k+1+4/n$.
Since
$$
\ell_kp_2\ge\left(1+\frac{4}{n}\right)p_2\ge p_*,
$$
for any $r\in[p_2,\infty]$, we see that
\begin{equation}
\label{eq:6.5}
\begin{split}
\|f(v(s))\|_{L^r}
&
\le \sum_{k=0}^{\infty} \frac 1{k!} \|v(s)\|_{L^{\ell_kr}}^{\ell_k}
\\
&
\le  \sum_{k=0}^{\infty} \frac 1{k!} \left( C s^{-\frac{n}{2}(\frac{1}{p_*} -\frac{1}{\ell_kr})} K\right )^{\ell_k}
\\
&
\le (CK)^{1+\frac{4}{n}}s^{\frac{n}{2r}-\frac{n}{2p_*}(1+\frac{4}{n})}
\sum_{k=0}^\infty\frac{\left(Cs^{-\frac{n}{2p_*}}K\right)^{2k}}{k!}
\\
&
\le (CK)^{1+\frac{4}{n}}s^{-\frac{n}{2}(\frac{1}{p_*}-\frac{1}{r})-\frac{2}{p_*}}\exp\left((CK)^2s^{-\frac{n}{p_*}}\right)
\end{split}
\end{equation}
for all $s>0$.
We can choose a sufficiently large constant $T_1\ge1$ such that, for all $s>T_1$,
it holds that
$$
\exp\left((CK)^2s^{-\frac{n}{p_*}}\right)\le 2.
$$
It is enough to chose
\be\label{eq:6.6}
T_1 \ge \left( \frac{(C K)^2}{\log 2}\right )^{\frac{p_*}{n}}.
\ee
This together with \eqref{eq:6.5} implies \eqref{eq:6.3}.
Thus Lemma~\ref{Lemma:6.1} follows.
$\Box$\vspace{7pt}

\noindent Let us prove now \eqref{eq:2.5} for small times.
\begin{lemma}
\label{Lemma:6.2}
Let $n\ge 1$ and $M = \|\varphi\|_{\exp L^2}<\eps$ small enough so that existence and uniqueness theorems apply.
Furthermore, let $u$ be the unique solution of the Cauchy problem \eqref{eq:1.1} satisfying \eqref{eq:1.8}.
Suppose $\varphi\in L^p$ for $p\in[1,2)$ and let $p_*$ be the constant given in \eqref{eq:2.3}.
Then, for any fixed $T_2\ge1$,
there exists a sufficiently small constant $\eps = \eps (p_*, T_2) >0$ such that,
if $M <\eps$,
then the solution $u$ satisfies
\begin{equation}
\label{eq:6.7}
\|u(t)\|_{L^q} \leq \tilde C t^{-\frac n2\left(\frac 1{p_*}-\frac 1q\right )}K,\quad 0<t\leq 2T_2
\end{equation}
for all $q\in[p_*,\infty]$ and for $\tilde C>0$ depending on $n$,
where $K$ is the constant defined by \eqref{eq:6.1}.
\end{lemma}
\noindent{\bf Proof.}
Let us consider
\[
u(t)= e^{t\Delta} \varphi + \int_0^t e^{(t-s)\Delta} f(u(s)) \, ds := e^{t\Delta} \varphi + F(u(t)).
\]
As for the linear part,
for any $q\in[p_*,\infty]$,
by \eqref{eq:3.3} and \eqref{eq:6.1} we have
\begin{equation}
\label{eq:6.8}
\|e^{t\Delta} \varphi \|_{L^q} \le c_2 t^{-\frac{n}{2}(\frac{1}{p_*} -\frac{1}{q})}\|\varphi \|_{L^{p_*}}
\le c_2 t^{-\frac{n}{2}(\frac{1}{p_*} -\frac{1}{q})}K,\quad t>0.
\end{equation}
Let us consider now  the nonlinear part $F(u)$.
For $r\ge p_1$ (defined in \eqref{eq:2.7}) such that $2r>n$ (and so $r$ depends on $n$ only),
we have
\[
\|F(u(t))\|_{L^\infty}
\le c_2 \int_0^{t} (t-s)^{-\frac n{2r}} \|f(u(s))\|_{L^r} \, ds
\le C t^{1-\frac n{2r}}\sup_{t>0}\|f(u(t))\|_{L^r},
\qquad t>0,
\]
where $C>0$ depends on $n$ only.
Now, by Lemma \ref{Lemma:3.4},
for $M<1$ small enough depending only on $r$,
we get
\[
C t^{1-\frac n{2r}} \sup_{t>0}\|f(u(t))\|_{L^r}
\le  C  t^{1-\frac n{2r}} r^3 M^{1+\frac{4}{n}}
\le  C T^{1-\frac n{2r}} M^{1+\frac{4}{n}},
\qquad t\leq 2T_2,
\]
where $C>0$ depends on $n$ only.
On the other hand,
since it follows from $p\in[1,2)$ with \eqref{eq:2.7} that
$$
1-\frac{n}{2}\left(\frac{1}{p_1}-\frac{1}{p_*}\right)>0,
$$
again by Lemma \ref{Lemma:3.4},
for $M<1$ small enough depending only on $p_1$,
we have
\begin{equation*}
\begin{split}
\|F(u(t))\|_{L^{p_*}}
&
\le c_2 \int_0^{t}(t-s)^{-\frac{n}{2}(\frac{1}{p_1}-\frac{1}{p_*})}\|f(u(s))\|_{L^{p_1}} \, ds
\\
&
\le C p_1^3t^{1-\frac{n}{2}(\frac{1}{p_1}-\frac{1}{p_*})} M^{1+\frac{4}{n}}
\le  C T_2^{1-\frac{n}{2}(\frac{1}{p_1}-\frac{1}{p_*})} M^{1+\frac{4}{n}},
\qquad t\leq 2T_2,
\end{split}
\end{equation*}
where $C>0$ depends on $n$.
If we choose
$$
M^{\frac{4}{n}} \max (T_2^{1-\frac{n}{2}(\frac{1}{p_1}-\frac{1}{p_*})}, {T_2}^{1-\frac{n}{2r}} )
<{T_2}^{-\frac{n}{2p_*}},
$$
then, due to \eqref{eq:6.2}, for any $q\in [p_*,\infty]$,  we get
 \begin{equation}
 \label{eq:6.9}
 \|F(u(t))\|_{L^q}
 \le  \tilde C MT_2^{-\frac{n}{2p_*}}\leq  \tilde C KT_2^{-\frac{n}{2p_*}},\qquad t\le 2T_2,
 \end{equation}
 where $\tilde C>0 $ depends on $n$.
Since $T_2\ge1$, gathering  \eqref{eq:6.8} and \eqref{eq:6.9},
for any $q\in [p_*,\infty]$,
we obtain
\[
\|u(t)\|_{L^q}
\le \tilde C \left( t^{-\frac n2\left(\frac 1p_*-\frac 1q\right )} + T_2^{-\frac{n}{2p_*}} \right )K
\le \tilde C \left( t^{-\frac n2\left(\frac 1p_*-\frac 1q\right )}
+ T_2^{-\frac{n}{2}(\frac{1}{p_*}-\frac{1}{q})} \right )K,
\quad t\le2T_2.
\]
Since $(2T_2)/t \geq 1$ we also get
\[
\|u(t)\|_{L^q} \leq  \tilde C t^{-\frac n2\left(\frac 1p_*-\frac 1q\right )}K,\qquad t\leq 2T_2
\]
for $\tilde C$ depending on $n$.
This implies \eqref{eq:6.7}.
$\Box$\vspace{7pt}
In order to prove estimate \eqref{eq:2.5} for all times $t>0$
we introduce once again a recursive sequence which, by uniqueness, converges to the solution $u$.
Let
\be\label{eq:6.10}
\left\{
\begin{aligned}
&u_0= e^{t\Delta}\varphi\\
&u_{j+1} = u_0 +  \int_0^t e^{(t-s)\Delta} f(u_j(s))\, ds,\quad j\geq 0.
\end{aligned}
\right .
\ee
Then, for $\varphi \ge0$,
 $u_j$ is an increasing sequence, namely,
$$
u_j(t,x)\le u_{j+1}(t,x)\le u(t,x),\qquad t>0,\quad x\in{\mathbb R}^n,
$$
for all $j\in{\mathbb N}_0$.
So, if $u$ is a solution of the Cauchy problem  \eqref{eq:1.1}
and satisfies \eqref{eq:1.8},
then, for all $j\in\N$, $u_j$ also satisfies
\be
\label{eq:6.11}
\sup _{t>0} \|u_j(t)\|_{\exp{L^2}} \leq 2\|\varphi\|_{\exp{L^2}}.
\ee
We end up by proving now Theorem \ref{Theorem:2.3}.

\bigskip
\noindent{\bf Proof of Theorem \ref{Theorem:2.3}.}
We are going to prove the first estimate \eqref{eq:2.5} uniformly for all $u_j(t,x)$
defined as in \eqref{eq:6.10} and so the estimate will pass to the limit $u(t,x)$.
Let $T_*$ be a sufficiently large constant to be chosen later, which satisfies
\begin{equation}
\label{eq:6.12}
T_*\ge T_1 \geq 1.
\end{equation}
Here $T_1$ is the constant given in Lemma~\ref{Lemma:6.1}, and depends only on $K$ and $p_*$.
First of all let us consider $u_0 =e^{t\Delta} \varphi$.
For any $q\in [p_*,\infty]$, we have
\be
\label{eq:6.13}
\|u_0(t)\|_{L^q} \leq c_2 t^{-\frac n2\left( \frac 1{p_*} -\frac 1q\right)}\|\varphi\|_{L^{p_*}}, \quad t>0.
\ee
By the monotonicity of $\{u_j\}$ and Lemma~\ref{Lemma:6.2} with $T_2=T_*$
we also have
\be
\label{eq:6.14}
\|u_j(t)\|_{L^q} \leq \tilde C t^{-\frac n2\left(\frac 1{p_*}-\frac 1q\right )}K,\quad 0<t\leq 2T_*,
\ee
for all $j\in\N$.
Let us assume that, for $j=m$,
\be\label{eq:6.15}
\|u_m(t)\|_q \leq 2 C_* t^{-\frac n2 \left( \frac 1{p_*} -\frac 1q\right )}K,\quad t>0,
\ee
with
$$
C_*=\max\{c_2,\tilde C\},
$$
and let us prove it for $j=m+1$.
Since
\[
u_{m+1}(t)= u_0(t) + \int_{0}^t e^{(t-s)\Delta} f(u_m(s))\,ds
\]
we consider only the nonlinear term.  For $q\in [p_*,\infty]$ and for  $t>2T_*$,
we split
\begin{equation}
\label{eq:6.16}
\begin{split}
&
\left\|
\int_{0}^t e^{(t-s)\Delta} f(u_m(s))\,ds \right \|_{L^q}
\\
&\le
 \int_{0}^{t/2}\left\| e^{(t-s)\Delta} f(u_m(s))\right \|_{L^q}\, ds
 + \int_{t/2}^{t}\left\| e^{(t-s)\Delta} f(u_m(s))\right \|_{L^q}\, ds
 \\
 & =: J_1(t) +J_2(t).
  \end{split}
\end{equation}
For $t>2T_*$ we have
\be
\label{eq:6.17}
J_1(t)
\le c_2
\left(
\int_{0}^{T_*}+\int_{T_*}^{t/2}\right)
(t-s)^{-\frac{n}{2}\left(\frac{1}{p_1}-\frac{1}{q}\right )}\|f(u_m(s))\|_{L^{p_1}}\,ds
:= A(t)+B(t),
\ee
where $p_1$ is the constant given in \eqref{eq:2.7}.
Due to estimate \eqref{eq:6.11} we can apply Lemma \ref{Lemma:3.4} to the $A(t)$ term,
and we obtain
\begin{equation}
\label{eq:6.18}
\begin{split}
A(t)
&
\le C t^{-\frac n2\left(\frac{1}{p_1}-\frac 1q\right )}\int_{0}^{T_*}M^{1+\frac{4}{n}} \, ds
\\
&
\le  Ct^{-\frac{n}{2}\left(\frac{1}{p_1}-\frac{1}{q}\right )} T_* M^{1+\frac{4}{n}}
=Ct^{-\frac{n}{2}\left(\frac{1}{p_*}-\frac{1}{q}\right )}t^{-\frac{n}{2}\left(\frac{1}{p_1}-\frac{1}{p_*}\right )}
 T_* M^{1+\frac{4}{n}}.
 \end{split}
 \end{equation}
As for the $B(t)$ term, we apply Lemma \ref{Lemma:6.1}.
Since it follows from \eqref{eq:2.3}, \eqref{eq:2.7} and \eqref{eq:6.4} with $p\in[1,2)$ that
\begin{equation*}
\left\{
\begin{array}{ll}
p_1=p_2=1
& \mbox{for}\quad 1\le n\le 4,
\vspace{5pt}\\
{\displaystyle
p_2=\max\left\{\frac{p_*n}{n+4},1\right\}<\frac{2n}{n+4}=p_1}
&\mbox{for}\quad n\ge5,
\end{array}
\right.
\end{equation*}
by \eqref{eq:6.3} we have
\begin{equation*}
 \begin{split}
 B(t)
 &
 \leq  C t^{-\frac n2\left(\frac{1}{p_1}-\frac 1q\right )}
 \int_{T_*}^{t/2}  K^{1+\frac{4}{n}} s^{-\frac n2\left(\frac 1{p_*} -\frac{1}{p_1}\right )-\frac{2}{p_*}}\,ds
 \\
 &
 \leq  C K^{1+\frac{4}{n}}
 t^{-\frac n2\left(\frac{1}{p_1}-\frac 1q\right )}
 \int_{T_*}^{t/2}s^{-\frac n2\left(\frac 1{p_*} -\frac{1}{p_1}\right )-\frac{2}{p_*}}\,ds.
\end{split}
\end{equation*}
For $p\in (p_1,2)$ (which implies $p_*=p$), we can choose $\sigma\in (0,1)$ satisfying
\be
\label{eq:6.19}
0<\sigma < \min\left\{\frac 2{p_*}-1,\frac n2\left(\frac{1}{p_1}-\frac 1{p_*}\right )\right\}.
\ee
So, we can write
\[
\begin{aligned}
B(t) &\leq
 C K^{1+\frac{4}{n}}
t^{-\frac n2\left(\frac 1{p_*}-\frac 1q\right )-\sigma}
\int_{T_*}^{t/2}   t^{-\frac n2\left(\frac{1}{p_1} -\frac1{p_*}\right )+\sigma }
s^{-\frac n2\left(\frac 1{p_*} -\frac{1}{p_1}\right )-\frac 2{p_*}}\,ds
\\
&
\le
 C K^{1+\frac{4}{n}}t^{-\frac n2\left(\frac 1{p_*}-\frac 1q\right )-\sigma}
\int_{T_*}^{t/2}  s^{-\frac n2\left(\frac{1}{p_1} -\frac1{p_*}\right )+\sigma }
s^{-\frac n2\left(\frac 1{p_*} -\frac{1}{p_1}\right )-\frac 2{p_*}}\,ds\\
&
\le  C K^{1+\frac{4}{n}}t^{-\frac n2\left(\frac 1{p_*}-\frac 1q\right )-\sigma}
\int_{T_*}^{\infty}   s^{-\frac 2{p_*}+\sigma }\,ds
\end{aligned}
\]
and in the end
$$
B(t)\le CK^{1+\frac{4}{n}}t^{-\frac n2\left(\frac 1{p_*}-\frac 1q\right )-\sigma}
T_*^{-\frac 2{p_*}+\sigma+1},
$$
where $C$ is a constant, independent of $M$, $K$ and $T_*$.
For  $p\le p_1$, namely $p_*=p_1$,
we cannot exploit the better decreasing in time as in \eqref{eq:6.19}.
Since $p_*<2$, for $p\le p_1$,
we get
\begin{equation}
\label{eq:6.20}
B(t) \leq  CK^{1+\frac{4}{n}}t^{-\frac n2\left(\frac{1}{p_*}-\frac 1q\right )}T_*^{-\frac{2}{p_*}+1}.
\end{equation}
In the end, for $p\in (p_1,2)$ we have the following estimates:
\[
J_1(t)\leq  C
\left( t^{-\frac n2\left(\frac{1}{p_1}-\frac 1q\right )} T_*M^{1+\frac{4}{n}}
+ K^{1+\frac{4}{n}}t^{-\frac n2\left(\frac 1{p_*}-\frac 1q\right )-\sigma}
T_*^{-\frac 2{p_*}+\sigma+1}
\right ),
\]
where $C$ is a constant, independent of $M$, $K$ and $T_*$.
Since $2/{p_*}-1-\sigma >0$,
we can take a sufficiently large constant $T_* \geq 1$ so that
\[
C K^{\frac{4}{n}}
T_*^{-\frac 2{p_*}+\sigma+1} \leq \frac 14 C_*
\]
which means
\be\label{eq:6.21}
 T_* \geq \left(\frac {4C K^{\frac{4}{n}}}{C_*}\right )^{\frac 1 {\frac 2{p_*}-1-\sigma}}.
 \ee
 This together with \eqref{eq:6.12} implies that
 $T_*$ depends on $K$ and $p$ but not on $M$.
Then we can also take a sufficiently small constant $M$ so that
\[
C T _*M^{\frac{4}{n}}\leq \frac 14 C_*
\]
and this means
\be\label{eq:6.22}
M \leq \left(\frac{4CT_*}{C_*}\right)^{-\frac{n}{4}}.
\ee
In the end,
since $T_*\ge1$,
for any $p\in(p_1,2)$,
by \eqref{eq:6.1}
we have
\begin{equation}
\label{eq:6.23}
J_1(t)\leq C_* \left(\frac 14 t^{-\frac n2\left(\frac{1}{p_1}-\frac 1q \right )}M
+\frac 14 t^{-\frac n2\left(\frac 1{p_*}-\frac 1q\right )-\sigma} K  \right )
\leq \frac{C_*}{2}t^{-\frac n2\left(\frac 1{p_*}-\frac 1q \right )}K
\end{equation}
for all $t>2T_*$.
In the $p\le p_1$ case we get from \eqref{eq:6.17}, \eqref{eq:6.18} and \eqref{eq:6.20}
similar conditions on $T_*$ and on $M$ in order to obtain  \eqref{eq:6.23}.

Let us come back to the $J_2(t)$ term in \eqref{eq:6.16}.
Since $q\ge p_*\ge p_2$,
by Lemma \ref{Lemma:6.1} again,
we have
\begin{equation*}
\begin{split}
J_2(t)
&
\leq  \int_{t/2}^{t}\left\|f(u_j(s))\right \|_{L^q}\, ds
\\
&
\leq   C  \int_{t/2}^{t} K^{1+\frac{4}{n}} s^{-\frac n2\left(\frac 1{p_*} -\frac 1q\right )-\frac 2{p_*}}\,ds
\leq CK^{1+\frac{4}{n}} t^{-\frac n2\left(\frac 1{p_*} -\frac 1q\right )-\frac{2}{p_*}+1}
\end{split}
\end{equation*}
Once again,
since $p_*<2$,
 we can choose $\lambda >0$ satisfying
$$
0<\lambda <\frac 2{p_*}-1,
$$
and we get
$$
J_2(t) \leq  CK^{1+\frac{4}{n}}
t^{-\frac n2\left(\frac 1{p_*} -\frac 1q\right )} t^{-\frac 2{p_*}+1+\lambda}
T_*^{-\lambda},\qquad t>2T_*,
$$
where $C$ is a constant, independent of $M$, $K$ and $T_*$.
For
\[
 CK^{\frac{4}{n}}T_*^{-\lambda} \leq \frac 14 C_*,
\]
which is implied for example by
\be
\label{eq:6.24}
T_*\geq \left(\frac {4C K^{\frac{4}{n}}}{C_*}\right)^{\frac 1\lambda},
\ee
we get for $t>2T_*$
\be\label{eq:6.25}
J_2(t) \leq   \frac {C_*}2 t^{-\frac n2\left(\frac 1{p_*} -\frac 1q\right )}   t^{-\frac 2{p_*} +1+\lambda }K.
\ee
This together with $T_*\ge1$ yields
\be\label{eq:6.26}
J_2(t) \leq  \frac{C_*}2 t^{-\frac n2\left(\frac 1{p_*} -\frac 1q\right )} K,\qquad t\ge2T_*.
\ee
Collecting \eqref{eq:6.1}, \eqref{eq:6.13}, \eqref{eq:6.23} and \eqref{eq:6.26} we have
\[
\begin{aligned}
\|u_{m+1}(t)\|_{L^q} &\leq c_2 t^{-\frac n2\left(\frac 1{p_*}-\frac 1q\right )} \|\varphi\|_{L^{p_*}}
+ C_* t^{-\frac n2\left(\frac 1{p_*} -\frac 1q\right )}K \\
&\leq 2 C_* t^{-\frac n2\left(\frac 1{p_*} -\frac 1q\right )}K,\qquad t>2T_*.
\end{aligned}
\]
This together with \eqref{eq:6.14} yields \eqref{eq:6.15} with $j=m+1$.
In order to make clear the dependence of the choice we made on $T_*$ and $M$, we collect
below all  the conditions \eqref{eq:6.12}, \eqref{eq:6.21}, \eqref{eq:6.22},  \eqref{eq:6.24}
\[
\begin{aligned}
&T_*\geq \max(1, T_1),\\
& T_* \geq \left(\frac {4CK^{\frac{4}{n}}}{C_*}\right )^{\frac 1 {\frac 2{p_*}-1-\sigma}},\\
&T_*\geq \left(\frac {4CK^{\frac{4}{n}}}{C_*}\right)^{\frac 1\lambda},\\
& M \leq \left(\frac {4CT_*}{C_*}\right)^{-\frac{n}{4}},
\end{aligned}
\]
with $T_1$ satisfying \eqref{eq:6.6}
\[
T_1 \ge \left( \frac{(C K)^2}{\log 2}\right )^{\frac{p_*}{n}},
\]
and $C_*$, $C$,  $\eta$, $\sigma$ and $\lambda$ constants depending at most on $n$ and $p_*$.
In the end, we find a function $F$ depending on $n$, $p_*$, $K$ such that the condition on $M$ can be written as
\[
M\leq F(n,p_*,K)
\]
as it was announced in the statement of Theorem \ref{Theorem:2.3}.
Finally we prove \eqref{eq:2.6}.
By \eqref{eq:6.23} and \eqref{eq:6.25}, for $p\in(p_1,2)$, we have
\[
t^{\frac n2\left(\frac 1{p_*}-\frac 1q\right ) } \|u(t)- e^{t\Delta} \varphi\|_{L^q} = o(1), \quad t\to \infty.
\]
Now, by density, let $\{\varphi_n\}\subset C_0^\infty $ such that $\varphi_n\to \varphi$ in $L^{p_*}$.
Then
\[
\begin{aligned}
t^{\frac n2\left(\frac 1{p_*}-\frac 1q\right ) } \|e^{t\Delta} \varphi\|_{L^q}
& \leq t^{\frac n2\left(\frac 1{p_*}-\frac 1q\right ) } \|e^{t\Delta} (\varphi-\varphi_n)\|_{L^q} +
t^{\frac n2\left(\frac 1{p_*}-\frac 1q\right ) } \|e^{t\Delta} \varphi_n\|_{L^q}\\
&\leq c_2 \left( \|\varphi-\varphi_n\|_{L^{p_*}}
+ t^{\frac n2\left(\frac 1{p_*}-\frac 1q\right ) } t^{-\frac n2\left(1-\frac 1q\right ) } \|\varphi_n\|_{L^1}\right )\\
&\leq  c_2 \left( \|\varphi-\varphi_n\|_{L^{p_*}}+ t^{-\frac n2\left(1-\frac 1{p_*}\right)} \|\varphi_n\|_{L^1}\right )
\end{aligned}
\]
for all $t>0$.
This proves  that
\[
t^{\frac n2\left(\frac 1{p_*}-\frac 1q\right ) } \| e^{t\Delta} \varphi\|_{L^q} = o(1), \quad t\to \infty
\]
and so
\[
t^{\frac n2\left(\frac 1{p_*}-\frac 1q\right ) } \|u(t)\|_{L^q} = o(1), \quad t\to \infty.
\]
Thus the proof of Theorem~\ref{Theorem:2.3} is complete.
$\Box$\vspace{7pt}
\begin{remark}\label{Remark:6.3}
It is worth commenting on the meaning of condition \eqref{eq:2.4}
appearing in Theorem~$\ref{Theorem:2.3}$ about the "smallness"
of the $\exp L^2$ norm of the initial data $\varphi$.
In fact, the evolution equation governing the Cauchy problem \eqref{eq:1.1} has no scaling invariance
and the $L^p$ and $\exp L^2$ norms have no relationship between each other.
In order to have initial data which fulfill  condition \eqref{eq:2.4},
let us choose a function $\varphi \in L^p \cap L^\infty$ with $p\in [1,2)$.
Since $\varphi\in L^2\cap L^\infty$, then $\varphi\in\exp L^2$ $($see {\rm\cite{A}}$)$.
Then, let us consider a dilation $\varphi_\lambda (x)= \lambda^{\frac np} \varphi(\lambda x)$
so that $\|\varphi_\lambda \|_{L^p}=\|\varphi\|_{L^p}$.
Since $\|\varphi_\lambda\|_{L^2}= \lambda^{n\left(\frac 1p -\frac 12\right )}\|\varphi\|_{L^2}$
and $\|\varphi_\lambda\|_{L^\infty}= \lambda^{\frac np}\|\varphi\|_{L^\infty}$,
it follows
\[
\limsup_{\lambda \to 0} \|\varphi_\lambda\|_{\exp L^2} 
\leq \lim_{\lambda \to 0}\left(\|\varphi_\lambda \|_{L^2} +\|\varphi_\lambda \|_{L^\infty}\right) =0.
\]
This implies that there is $\lambda >0$ so that $\varphi_\lambda$ fulfills condition \eqref{eq:2.4}, even though its $L^p$ norm might be large.
\end{remark}
\vspace{8pt}

We end  this section by proving Theorem~\ref{Theorem:2.4}.
In the following Lemma we assume $\|v(s)\|_{L^q}$ bounded at the origin and decaying at infinity,
and we can deduce that also  $\|f(v(s))\|_{L^r}$ is bounded and decays at infinity
for $r\geq p_2$ where $p_2$ is defined in \eqref{eq:6.4}.
\begin{lemma}
\label{Lemma:6.3}
Let $n\geq 1$, $p\in [1,2)$ and $L>0$.
Let  $v(t,x)$ a function satisfying
\begin{equation}
\label{eq:6.27}
\sup_{s>0} (1+s)^{\frac{n}{2}\left(\frac{1}{p}-\frac{1}{q}\right )}\|v(s)\|_{L^q} \leq C L
 \end{equation}
 for all $q\in[p,\infty]$ and $C$  independent of $q$.
Then,  there is $\delta >0$  such that  if $\tilde K<\delta$, then for all $r\in [p_2,\infty]$
\begin{equation}
\label{eq:6.28}
\sup_{s>0}  (1+s)^{\frac{n}{2}\left(\frac{1}{p}-\frac{1}{r}\right )+\frac{2}{p}}
\|f(v(s))\|_{L^r} \leq 2 (CL)^{1+\frac{4}{n}}\end{equation}
where $f$ and $p_2$ are defined as in \eqref{eq:1.6} and \eqref{eq:6.4}, respectively.
\end{lemma}
\noindent{\bf Proof.}
Let $r_1=4/n +1$.
Then,
since $r_1p_2\ge p$,
for any $r\in[p_2,\infty]$,
it follows from \eqref{eq:1.6}, \eqref{eq:3.4} and \eqref{eq:6.27} that
\begin{equation}
\label{eq:6.29}
\begin{split}
\|f(v(s))\|_{L^r}
&
\le
\sum_{k=0}^{\infty} \frac {\|v^{2k+r_1}(s)\|_{L^r}}{k!}
\\
&
\le
\sum_{k=0}^{\infty} \frac {\|v(s)\|_{L^{(2k+r_1)r}}^{2k+r_1}}{k!}
\\
&
\le
\sum_{k=0}^{\infty}
\frac {1}{k!}\left(C(1+s)^{-\frac{n}{2}\left(\frac{1}{p}-\frac{1}{(2k+r_1)r}\right)}L\right)^{2k+r_1}
\\
&
\le
(C\tilde K)^{r_1}(1+s)^{\frac{n}{2r}-\frac{nr_1}{2p}}\sum_{k=0}^\infty\frac{1}{k!}
\left(C(1+s)^{-\frac{n}{2p}}L\right)^{2k}
\\
&
\le
(C\tilde K)^{r_1}(1+s)^{\frac{n}{2r}-\frac{nr_1}{2p}}\sum_{k=0}^\infty\frac{1}{k!}\left(CL\right)^{2k}
\end{split}
\end{equation}
for all $s>0$.
We can take a sufficiently small $\delta$, which independent of $r$, so that,
for $L\leq \delta $,
it holds that
$$
\sum_{k=0}^\infty\frac{1}{k!}\left(CL\right)^{2k}
=e^{(C\tilde K)^2}\le2.
$$
This together with \eqref{eq:6.29} implies \eqref{eq:6.28}
Thus Lemma~\ref{Lemma:6.3} follows.
$\Box$\vspace{7pt}
\noindent{\bf Proof of Theorem \ref{Theorem:2.4}.}
Put $L=\max\{\|\varphi\|_{\exp L^2},\|\varphi\|_{L^p}\}$.
Applying the same argument as in the proof of Theorem~\ref{Theorem:2.2} with
Lemma~\ref{Lemma:6.3},
we can prove Theorem~\ref{Theorem:2.4}.
So we omit the proof.
$\Box$\vspace{7pt}
%

\section{Asymptotic behavior}

Let us come to the asymptotic behavior of the solution $u$ as stated in Theorem \ref{Theorem:2.5}.

\bigskip
\noindent{\bf Proof of Theorem \ref{Theorem:2.5}.}
By the assumption for the initial data $\varphi$ with \eqref{eq:1.8}
we can apply Lemma~\ref{Lemma:3.4} and obtain
$$
\sup_{t>0}\|f(u(t))\|_{L^1}<\infty,\qquad 1\le n\le 4.
$$
For the case $n\ge5$, since
$$
f(u(t))=u^{1+\frac{4}{n}}(t)\,e^{u^2(t)}=u^{1+\frac{4}{n}}(t) +  u^{1+\frac{4}{n}}(t)\left(e^{u^2(t)}-1\right ),
$$
by a proof analogous to Lemma \ref{Lemma:3.4} and \eqref{eq:2.9} we get
 $$
\sup_{t>0}\|f(u(t))\|_{L^1}\le
\sup_{t>0}\left(\|u(t)\|_{L^{1+\frac{4}{n}}}^{1+\frac{4}{n}}
+\left\|u^{1+\frac{4}{n}}(t)\left(e^{u^2(t)}-1\right )\right\|_{L^1}\right)<\infty.
$$
Therefore we can define a mass of $u(t)$ denote by $m(t)$, that is,
$$
m(t):=\int_{{\mathbb R}^n}\varphi(x)\,dx+\int_0^t\int_{{\mathbb R}^n}f(u(s,x))\,dx\,ds,\qquad t\ge0.
$$
In fact, if $\varphi\in L^1$, then we can define the integral equation \eqref{eq:1.7} in $L^1$
for almost all $t>0$,
and by ($G1$) we have
\begin{equation*}
\begin{split}
&
\int_{\R^n} e^{t\Delta} \varphi (x)\,dx = \int_{\R^n} \varphi (x)\,dx,
\\
&
\int_{\R^n} \int_0^t \e^{(t-s)\Delta} f(u(s,x))\,ds\,dx = \int_0^t\int_{{\mathbb R}^n}  f(u(s,x))\,dx\,ds.
\end{split}
\end{equation*}
Let $T_1$ be the constant given in Lemma~\ref{Lemma:6.1}.
Then, by Lemma~\ref{Lemma:6.1} we have
\begin{equation}
\label{eq:7.1}
\begin{split}
 \int_0^t \int_{\R^n}  f(u(s,x))\,dx\,ds
&
=  \int_0^{T_1}\int_{\R^n} f(u(s,x))\,dx\,ds
+  \int_{T_1}^t\int_{\R^n} f(u(s,x))\,dx\, ds
\\
&
\le CT_1 + C \int_{T_1}^t s^{-2}\,ds
\le C(T_1+T_1^{-1})
\end{split}
\end{equation}
for all $t\ge T_1$.
This implies that there exists the limit of $m(t)$, which we denote by $m_*$, such that
\[
m_*:=
\lim_{t\to \infty}m(t)
= \int_{\R^n} \varphi (x) dx +\int_0^{\infty}\int_{{\mathbb R}^n}  f(u(s,x))\,dx\,ds.
\]
Furthermore, similarly to \eqref{eq:7.1}, we obtain
\[
m_*-m(t)
\le C\int_t^\infty\int_{\R^n}f(u(s,x))\,dx\,ds
\le Ct^{-1}
\]
for all $t\ge T_1$.
Therefore, applying an argument similar to the proof of Theorem~1.1 in \cite{FIK2012}
with \eqref{eq:2.8} and \eqref{eq:2.9}, we have \eqref{eq:2.10},
and the proof of Theorem~\ref{Theorem:2.5} is complete.

\hfill $\Box$\vspace{7pt}
%

\section{Generalization}
In this last section, we are going to get similar results as in the Section~2 to the general problem.
Let $\theta\in(0,2]$ and $r>1$.
Put
\begin{equation}
\label{eq:8.1}
f(u)=|u|^{\frac{r\theta}{n}}ue^{u^r}.
\end{equation}
Then, similarly to Lemma~\ref{Lemma:3.4}, the following hold.
\begin{lemma}
\label{Lemma:8.1}
Let $n\geq 1$ and $M>0$.
Assume that a function $u\in L^\infty(0,\infty\,;\, \exp L^r)$ satisfies the condition
$$
\sup _{t>0} \|u(t)\|_{\exp{L^r}} \leq M.
$$
Let $f$ be the function defined as in \eqref{eq:8.1}.
Put
$$
p_1=\max\{rn/(n+r\theta),1\}.
$$
Then, for all $p\in [p_1,\infty)$ there is $\eps =\eps (p) >0$ such that if $M<\eps$, then
$$
\sup_{t>0} \|f(u(t))\|_{L^p} \leq 2Cp^3M^{1+\frac{r\theta}{n}},
$$
where $C$ is independent of $p$, $n$ and $M$.
\end{lemma}
Let us consider now for simplicity only the integral equation
\be\label{eq:8.2}
u(t)= e^{-t{\cal L}_\theta} \varphi + \int_0^t e^{-(t-s) {\cal L}_\theta} f(u(s)) \, ds.
\ee
This is the integral formulation of \eqref{eq:2.11}.
In a similar way as in Section~4 
one can prove that this integral equation is equivalent to the differential equation 
if we consider small solutions.

Applying the arguments as in the previous sections with Lemma~\ref{Lemma:8.1},
we have the following.
\begin{theorem}
\label{Theorem:8.1}
Let $n\ge1$ and $r>1$.
Assume $\varphi \in \exp{L^r}$.
Then there exists $\varepsilon =\eps(n) >0$  such that,
if $\|\varphi\|_{\exp{L^r}} <\eps$,
then there exists a unique solution $u$ of the integral equation \eqref{eq:8.2} satisfying
$$
u \in L^\infty(0,\infty\,;\, \exp{L^r})
$$
and
\begin{equation}
\label{eq:8.3}
\sup _{t>0} \|u(t)\|_{\exp{L^r}} \leq 2\|\varphi\|_{\exp{L^r}}.
\end{equation}
\end{theorem}
\begin{theorem} 
\label{Theorem:8.2}
Let $n\ge1$, $r>1$ and $\varphi\in\exp L^r$ with $\varphi \geq 0$.
Assume that there exists a unique positive solution $u$ of \eqref{eq:8.2} satisfying \eqref{eq:8.3}. 
Then there exist  $\eps = \eps (n) >0$ and $C=C(n)>0$ such that, 
if $\|\varphi\|_{\exp L^r} <\eps$, 
then the solution $u$ satisfies
\[
\|u(t)\|_{L^q} \leq  C t^{-\frac n\theta \left( \frac 1r -\frac 1q\right )}\|\varphi\|_{\exp L^r},\quad t>0,
\]
for any $q\in[r,\infty]$.
\end{theorem}
\begin{theorem} 
\label{Theorem:8.3}
Assume the same conditions as in Theorem~$\ref{Theorem:8.2}$.
Furthermore, suppose that $\varphi\in L^p$ for some $p\in[1,r)$.
Put
\[
p^*=\max\left\{p,\frac{rn}{n+r\theta}\right\}.
\]
Then there exist positive constants $\eps =\eps (n) $, $C=C(n)$
and a positive function $F=F(n, p^*, \|\varphi\|_{L^{p^*}})$ such that, if
\[
\|\varphi\|_{\exp L^r} < \min \left(\eps, F(n,p^*, \|\varphi\|_{L^{p^*}})\right),
\]
then the solution $u$ satisfies
\be\label{eq:8.4}
\|u(t)\|_{L^q} 
\leq C t^{-\frac n\theta \left( \frac 1{p^*} -\frac 1q\right )}\|\varphi\|_{\exp L^r\cap L^{p^*}},\quad t>0,
\ee
for any $q\in[p^*,\infty]$.
In particular, if $p\in (p_3,r)$, then
\be\label{eq:8.5}
\|u(t)\|_{L^q} = o\left(t^{-\frac n\theta \left( \frac 1{p^*} -\frac 1q\right )}\right ),\quad t\to \infty.
\ee
Here $[K]$ is the integer satisfying $K-1\le[K]< K$ and
$$
p_3:=\max\left\{1,\frac{rn}{n+r\theta}\right\}.
$$
\end{theorem}
\begin{theorem} 
\label{Theorem:8.4}
Assume the same conditions as in Theorem~$\ref{Theorem:8.3}$.
Then there exists a positive constant $\delta=\delta(n)$ such that,
if
\[
\max\{\|\varphi\|_{\exp L^r},\|\varphi\|_{L^p}\}<\delta,
\]
then \eqref{eq:8.4} with $p^*=p$ holds for all $q\in[p,\infty]$
In particular, for all $q\in[p,\infty)$,
\[
\|u(t)\|_{L^q} \leq C(1+t)^{-\frac n\theta \left( \frac 1p -\frac 1q\right )}\|\varphi\|_{\exp L^r\cap L^p},\quad t>0.
\]
Furthermore, if $p\in (1,r)$, then \eqref{eq:8.5} with $p^*=p$ holds.
\end{theorem}
\begin{theorem} 
\label{Theorem:8.5}
Let $n\ge 1$, $\varphi \ge 0$ and $\varphi \in \exp L^r \cap L^1$.
Assume $\|\varphi\|_{\exp L^r}$ is small enough.
Furthermore, suppose that
\begin{enumerate}[a)]
\item for $n\geq 1$,
\[
\sup_{t>0}\,t^{\frac{n}{\theta}\left(1 -\frac 1q\right )}\|u(t)\|_{L^q}<\infty,
\qquad q\in[1,\infty];
\]
\item for $n>\max\{1,[r\theta/(r-1)]\}$, assume moreover that there is $T^*>0$ such that
\[
\sup_{0<t<T^*}\,\|u(t)\|_{L^{\frac{r\theta}{n}+1}}<\infty.
\]
\end{enumerate}
Then there exists the limit
\[
 \lim_{t\to \infty} \int_{\R^n}u(x,t) \,dx
 = \int_{\R^n} \varphi(x)\,dx + \int_{0}^{\infty} \int_{\R^n} f(u(t,x))\,dx\,dt :=m_*
\]
such that
 \[
 \lim_{t\to \infty} t^{\frac{n}{\theta}\left(1 -\frac 1q\right )}\|u(t)-m_* G_\theta(t+1)\|_{L^q}=0,
 \]
 for any $q\in[1,\infty]$.
\end{theorem}

\noindent
{\bf Acknowledgements.}
The work of the second author (T. Kawakami) was supported in part 
by Grant-in-Aid for Young Scientists (B) (No. 24740107) and (No. 16K17629) 
of JSPS (Japan Society for the Promotion of Science)
and by the JSPS Program for Advancing Strategic International Networks 
to Accelerate the Circulation of Talented Researchers 
``Mathematical Science of Symmetry, Topology and Moduli, 
Evolution of International Research Network based on OCAMI''.



\bibliographystyle{amsplain}

\par \vspace{1cm} \noindent
{\it Adresses:}
\par \medskip \noindent
\author{Giulia Furioli \\
{\small DIGIP, Universit\`a
di Bergamo,}\\
{\small Viale Marconi 5, I--24044 Dalmine (BG), Italy} \\
{\small E-mail:} \texttt{gfurioli@unibg.it}\\
\ \\
Tatsuki Kawakami\\
{\small Department of Mathematical Sciences,}\\
{\small Osaka Prefecture University,}\\
{\small Gakuencho 1-1, Sakai 599-8531, Japan} \\
{\small E-mail: } \texttt{kawakami@ms.osakafu-u.ac.jp}\\
\
\\
 Bernhard Ruf\\
{\small Dipartimento di Matematica F.~Enriques,}\\ 
{\small Universit\`a degli studi di Milano,}\\
{\small Via Saldini 50 , I--20133 Milano, Italy} \\
{\small E-mail:} \texttt{bernhard.ruf@unimi.it}\\
\
\\
 Elide Terraneo\\
{\small Dipartimento di Matematica F.~Enriques,}\\ 
{\small Universit\`a degli studi di Milano,}\\
{\small Via Saldini 50 , I--20133 Milano, Italy} \\
{\small E-mail:} \texttt{elide.terraneo@unimi.it}\\
}
\end{document}